\newif\iffinal
\finalfalse	% Not final version
\finaltrue	% Final version

\documentclass[letterpaper,12pt,reqno]{amsart}

\usepackage{mathtools}
\mathtoolsset{showonlyrefs=true}
\mathtoolsset{showonlyrefs=false}

\RequirePackage[utf8]{inputenc}
\usepackage[portrait,margin=3cm]{geometry}
\iffinal\else\usepackage[notref,notcite]{showkeys}\fi
\usepackage{mathrsfs,soul}
\usepackage[hyperfootnotes=false]{hyperref}
\usepackage{harpoon}

\usepackage[foot]{amsaddr}
\usepackage{amssymb,amsmath,amsthm,amsfonts,amsbsy,latexsym,dsfont}
\usepackage{leftidx}
\usepackage[textsize=small]{todonotes}       % includes TO DO LIST
\usepackage{graphicx}
\usepackage[numeric,initials,nobysame]{amsrefs}
\usepackage{upref,setspace}
\usepackage{xcolor,colortbl}
\usepackage{enumerate}
\newenvironment{enumeratei}{\begin{enumerate}[\upshape (i)]}{\end{enumerate}}
\newenvironment{enumeratea}{\begin{enumerate}[\upshape (a)]}{\end{enumerate}}

\usepackage{tikz}
\usetikzlibrary{calc,intersections,through,backgrounds,shapes.geometric}
\usetikzlibrary{graphs}
 	\usepackage{pgfplots}

\makeatletter
\DeclareFontFamily{U}{MnSymbolA}{}
\DeclareFontShape{U}{MnSymbolA}{m}{n}{
    <-6>  MnSymbolA5
   <6-7>  MnSymbolA6
   <7-8>  MnSymbolA7
   <8-9>  MnSymbolA8
   <9-10> MnSymbolA9
  <10-12> MnSymbolA10
  <12->   MnSymbolA12}{}
\DeclareFontShape{U}{MnSymbolA}{b}{n}{
    <-6>  MnSymbolA-Bold5
   <6-7>  MnSymbolA-Bold6
   <7-8>  MnSymbolA-Bold7
   <8-9>  MnSymbolA-Bold8
   <9-10> MnSymbolA-Bold9
  <10-12> MnSymbolA-Bold10
  <12->   MnSymbolA-Bold12}{}
\DeclareSymbolFont{MnSyA}{U}{MnSymbolA}{m}{n}
\SetSymbolFont{MnSyA}{bold}{U}{MnSymbolA}{b}{n}

\DeclareRobustCommand{\overleftharpoon}{\mathpalette{\overarrow@\leftharpoonfill@}}
\DeclareRobustCommand{\overrightharpoon}{\mathpalette{\overarrow@\rightharpoonfill@}}
\def\leftharpoonfill@{\arrowfill@\leftharpoondown\mn@relbar\mn@relbar}
\def\rightharpoonfill@{\arrowfill@\mn@relbar\mn@relbar\rightharpoonup}

\DeclareMathSymbol{\leftharpoondown}{\mathrel}{MnSyA}{'112}
\DeclareMathSymbol{\rightharpoonup}{\mathrel}{MnSyA}{'100}
\DeclareMathSymbol{\mn@relbar}{\mathrel}{MnSyA}{'320}
\makeatother

\usepackage{paralist}

\newenvironment{inparaenumi}{\begin{inparaenum}[\upshape \bfseries (i) ]}{\end{inparaenum}}

\newenvironment{inparaenuma}{\begin{inparaenum}[\upshape \bfseries (a) ]}{\end{inparaenum}}

\numberwithin{equation}{section}
\numberwithin{figure}{section}
\numberwithin{table}{section}

\usepackage{caption}
\usepackage{subcaption}

\newtheorem{thm}{Theorem}[section]
\newtheorem{lem}[thm]{Lemma}

\newtheorem{cor}[thm]{Corollary}

\newtheorem{prop}[thm]{Proposition}

\newtheorem{conj}[thm]{Conjecture}

\newtheorem{constr}[thm]{Construction}

\theoremstyle{definition}
\newtheorem{rem}{Remark}

\renewcommand{\leq}{\le}
\renewcommand{\geq}{\ge}

\newcommand{\ind}{\mathds{1}}
\newcommand{\eps}{\varepsilon}

\newcommand{\equald}{\stackrel{\mathrm{d}}{=}}

\newcommand{\weakc}{\stackrel{\mathrm{d}}{\longrightarrow}}
\newcommand{\convas}{\stackrel{\mathrm{a.s.}}{\longrightarrow}}

\newcommand{\tonn}{\stackrel{n\to\infty}{\xrightarrow{\hspace*{.8cm}}}}

\newcommand{\hght}{\mathrm{ht}}
\newcommand{\bm}{\mathrm{bm}}
\newcommand{\br}{\mathrm{br}}
\newcommand{\exx}{\mathrm{ex}}

\newcommand{\radi}{\mathrm{Rad}}

\newcommand{\spls}{\mathrm{sp}}
\newcommand{\INT}{\mathrm{Int}}

\newcommand{\df}{\mathrm{DF}}
\newcommand{\dft}{\mathrm{DFT}}
\DeclareMathOperator{\bbf}{BF}
\DeclareMathOperator{\sbf}{bf}
\DeclareMathOperator{\bfac}{BFAC}
\DeclareMathOperator{\dfac}{DFAC}
\DeclareMathOperator{\bft}{BFT}

\newcommand{\crumg}{\mathrm{CRUM}_{(g)}}
\newcommand{\crumgbf}{\mathrm{CRUM}_{(g)}^{\mathrm{BF}}}
\newcommand{\bum}{\mathbb{UM}_{n, g}}
\newcommand{\bumast}{\mathbb{UM}_{n, g}^{\ast}}
\newcommand{\cur}{\mathrm{cur}}
\newcommand{\sort}{\mathrm{sort}}

\makeatletter
\DeclareFontFamily{U}{tipa}{}
\DeclareFontShape{U}{tipa}{m}{n}{<->tipa10}{}
\newcommand{\arc@char}{{\usefont{U}{tipa}{m}{n}\symbol{62}}}%

\newcommand{\arc}[1]{\mathpalette\arc@arc{#1}}

\newcommand{\arc@arc}[2]{%
  \sbox0{$\m@th#1#2$}%
  \vbox{
    \hbox{\resizebox{\wd0}{\height}{\arc@char}}
    \nointerlineskip
    \box0
  }%
}
\makeatother

\def\qed{ \hfill $\blacksquare$}
\newcommand{\cH}{\mathcal{H}}\newcommand{\cI}{\mathcal{I}}
\newcommand{\cL}{\mathcal{L}}

\newcommand{\cR}{\mathcal{R}}
\newcommand{\cT}{\mathcal{T}}

\newcommand{\ve}{\mathbf{e}}

\newcommand{\vm}{\mathbf{m}}

\newcommand{\vt}{\mathbf{t}}

\newcommand{\mvH}{\boldsymbol{H}}\newcommand{\mvI}{\boldsymbol{I}}

\newcommand{\mvX}{\boldsymbol{X}}

\newcommand{\mve}{\boldsymbol{e}}\newcommand{\mvf}{\boldsymbol{f}}
\newcommand{\mvh}{\boldsymbol{h}}\newcommand{\mvj}{\boldsymbol{j}}

\newcommand{\mvy}{\boldsymbol{y}}

\newcommand{\mvxi}{\boldsymbol{\xi}}

\newcommand{\fB}{\mathfrak{B}}\newcommand{\fC}{\mathfrak{C}}
\newcommand{\fD}{\mathfrak{D}}

\newcommand{\fN}{\mathfrak{N}}
\newcommand{\fR}{\mathfrak{R}}

\newcommand{\fX}{\mathfrak{X}}

\newcommand{\bD}{\mathbb{D}}\newcommand{\bE}{\mathbb{E}}
\newcommand{\bH}{\mathbb{H}}

\newcommand{\bM}{\mathbb{M}}
\newcommand{\bR}{\mathbb{R}}
\newcommand{\bS}{\mathbb{S}}\newcommand{\bU}{\mathbb{U}}

\newcommand{\bZ}{\mathbb{Z}}
\DeclareMathOperator{\E}{\mathbb{E}}
\DeclareMathOperator{\pr}{\mathbb{P}}

\DeclareMathOperator{\var}{Var}

\DeclareMathOperator{\dis}{dis}

\DeclareMathOperator{\GHP}{GHP}
 \DeclareMathOperator{\height}{ht}

\newcommand{\erdos}{Erd\H{o}s-R\'enyi }

\newcommand{\esbf}{\mathbf{e}_{(s)}^{\mathrm{BF}}}
\newcommand{\nsbf}{\eta_{(s)}^{\mathrm{BF}}}
\newcommand{\esdf}{\mathbf{e}_{(s)}^{\mathrm{DF}}}
\newcommand{\egum}{\mathbf{e}_{(g)}^{\mathrm{UM}}}

\newcommand{\hsbf}{\mathcal{H}_{(s)}^{\mathrm{BF}}}
\newcommand{\htwogbf}{\mathcal{H}_{(2g)}^{\mathrm{BF}}}

\newcommand{\degr}{\mathrm{deg}}

\usepackage{fourier}

\definecolor{aqua}{rgb}{0.0, 1.0, 1.0}
\definecolor{webbrown}{rgb}{.6,0,0}
\definecolor{pinegreen}{rgb}{0.0, 0.47, 0.44}
\definecolor{ultramarineblue}{rgb}{0.25, 0.4, 0.96}
\definecolor{jrnl}{rgb}{0.0, 0.5, 0.0}
\definecolor{lincolngreen}{rgb}{0.11, 0.35, 0.02}
\definecolor{green(html/cssgreen)}{rgb}{0.0, 0.5, 0.0}
\definecolor{airforceblue}{rgb}{0.36, 0.54, 0.66}
\definecolor{azure}{rgb}{0.0, 0.5, 1.0}
\definecolor{bleudefrance}{rgb}{0.19, 0.55, 0.91}
\definecolor{cobalt}{rgb}{0.0, 0.28, 0.67}

\hypersetup{colorlinks=true, linktocpage=true, pdfstartpage=1, pdfstartview=FitV,
    breaklinks=true, pdfpagemode=UseNone, pageanchor=true, pdfpagemode=UseOutlines,
    plainpages=false, bookmarksnumbered, bookmarksopen=true, bookmarksopenlevel=1,
    hypertexnames=true, pdfhighlight=/O,
    urlcolor=webbrown, linkcolor=cobalt, citecolor=jrnl
    }

\newcommand{\greg}[1]{\textcolor{black}{#1}}

\newcommand{\sanch}[1]{\textcolor{black}{{#1}}}

\newcommand{\chng}[1]{\textcolor{black}{#1}}

\begin{document}

\title[Breadth-first constructions]{On breadth-first constructions of scaling limits of random graphs and random unicellular maps}

\date{}
% Secondary: ;}

\author[Miermont]{Gr{\' e}gory Miermont$^{1}$}
\address{$^{1}$Unit\'{e} de Math\'{e}matiques Pures et Appliqu\'{e}es, \'{E}cole Normale Sup\'{e}rieure de Lyon, France}
\author[Sen]{Sanchayan Sen$^2$}
\address{$^2$Department of Mathematics, Indian Institute of Science, Bangalore, India}
\email{gregory.miermont@ens-lyon.fr, sanchayan.sen1@gmail.com}

\maketitle
\begin{abstract}
We give alternate constructions of (i) the scaling limit of the uniform connected graphs with given fixed surplus, and (ii) the continuum random unicellular map (CRUM) of a given genus
that start with a suitably tilted Brownian continuum random tree and make `horizontal' point identifications, at random heights, using the local time measures.
Consequently, this can be seen as a continuum analogue of the breadth-first construction of a finite connected graph.
In particular, this yields a breadth-first construction of the scaling limit of the critical \erdos random graph which answers a question posed in \cite{BBG-12}.
As a consequence of this breadth-first construction we obtain descriptions of the radii, the distance profiles, and the two point functions of these spaces in terms of functionals of tilted Brownian excursions.
\end{abstract}

\section{Introduction and definitions}
\label{sec:intro}

\greg{This paper studies properties of random metric spaces
  that arise naturally in the study of critical random graphs and
  random maps, by providing new constructions of these objects. 
  A common feature of these spaces is that they look
  locally like random trees, which makes it possible to construct them
by performing certain gluings in models of continuum random trees. One
nice aspect of the particular ``breadth-first construction'' considered
in this paper is that it can be defined canonically for a wide family
of (deterministic) metric spaces called $\bR$-graphs. We start by explaining this construction. }

\greg{Recall that a metric space $(X,d)$ is called geodesic if for every
$x,y\in X$, there exists an isometric embedding $\psi:[0,d(x,y)]\to X$
with $\psi(0)=x$ and $\psi(d(x,y))=y$. We call $\psi$ a geodesic path between $x$ and $y$. }
A compact geodesic metric space is called \greg{an $\bR$-tree}
 \cite{legall-survey,evans-book} if it has no embedded cycles.
A compact geodesic metric space $(X,d)$ is called an $\bR$-graph \cite{AddBroGolMie13} if for every $x\in X$, there exists $\eps>0$ such that 
\[
B(x,\eps; X):=\big\{y\in X\ :\ d(y, x)\leq\eps\big\}
\]
with the induced metric $d|_{B(x,\eps; X)}$ is \greg{an $\bR$-tree.}
A measured $\bR$-graph is an $\bR$-graph with a probability measure on its Borel $\sigma$-algebra.
By \cite[Theorem 2.7]{AddBroGolMie13}, $(X, d)$ is an $\bR$-graph if
and only if there exists a finite connected multigraph $G=(V, E)$ and
a collection $\{(T_e, x_e, y_e)\ :\ e\in E\}$ where $T_e$ is \greg{an $\bR$-tree} 
and $x_e, y_e\in T_e$ such that $(X, d)$ is isometric to the space constructed from $G$ by \greg{performing a metric gluing of the spaces $T_e$,} replacing each edge $e\in E$ by $T_e$ \greg{and} identifying $x_e$ with one endpoint of $e$ and $y_e$ with the other endpoint of $e$.
Random measured $\bR$-graphs arise naturally as scaling limits of various random graphs \cite{BBG-12, bhamidi-broutin-sen-wang, bhamidi-hofstad-sen, bhamidi-sen, SBSSXW14, lab-oa-gc-ef-cg}.

For an $\bR$-graph $(X, d)$ and $x\in X$, choose $\eps>0$ such that $B(x,\eps; X)$ is \greg{an $\bR$-tree}, and define the degree of $x$ as 
\[
\degr(x; X):=\#\big\{\text{connected components of }B(x,\eps; X)\setminus\{x\}\big\}.
\]
Note that the value of $\degr(x; X)$ is independent of the choice of $\eps$.
A triple $(X, d, x_\ast)$ where $(X, d)$ is an $\bR$-graph  and $x_\ast\in X$ is a distinguished point is called a rooted $\bR$-graph.
We can similarly define a measured rooted $\bR$-graph.
We define the radius of $X$ to be 
\[
\radi(X):=\sup_{x\in X}d(x, x_{\ast}).
\]

\begin{figure}
	\centering
   	\includegraphics[trim=0cm 0cm 0cm 0cm, clip=true, angle=0, scale=.31]{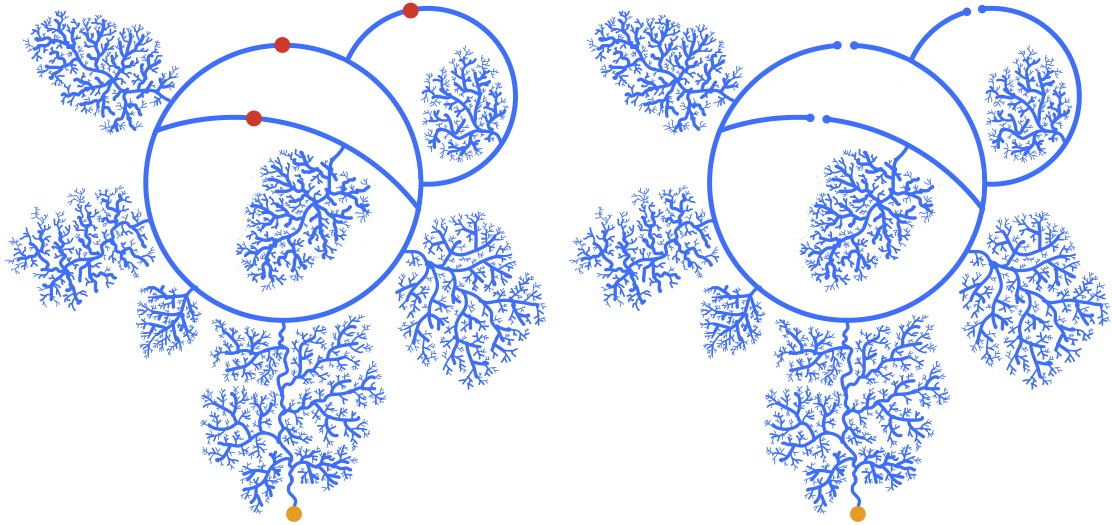}
	\captionof{figure}{On the left, an $\bR$-graph $X$ with its root colored orange and points in $\fX$ colored red. On the right, the space $X'$.}
	\label{fig:1-2}
\end{figure}

Let $(X, d, x_\ast)$ be a rooted $\bR$-graph.
Let $\fX$ be the cut locus of $x_\ast$ in X, i.e.,
the set of all points $x$ in $X$ such that there exist at least two distinct (not necessarily disjoint) geodesic paths from $x_\ast$ to $x$.
Then $X\setminus\fX$ is connected\greg{, since no point on the unique geodesic between a point $x\in X\setminus \fX$ and $x_*$ can be in $\fX$}.
Assume further that 
\begin{align}\label{eqn:1}
\degr(x; X)=2\ \text{ for every }\ x\in\fX.
\end{align}
Then it is easy to see that $\fX$ is a finite set\greg{, for instance from the description of $X$ as a metric gluing of doubly-rooted $\bR$-trees $(T_e,x_e,y_e)$, $e\in E$, along a finite graph $G=(V, E)$: for every point $x$ of $\fX$, there exists $e\in E$ such that $e$ belongs to one of the finitely many cycles in $G$ and $x$ is a point on the geodesic path connecting $x_e$ and $y_e$ in $T_e$, and removing $x$ has the effect of decreasing the total number of cycles.}
From now on we will only work with $\bR$-graphs that satisfy \eqref{eqn:1}.

For a rooted $\bR$-graph $(X, d, x_\ast)$ satisfying \eqref{eqn:1}, define $X'$ to be the completion of the metric space $X\setminus\fX$ \greg{endowed with the intrinsic metric inherited from $d$.}
Let $d'$ denote the metric on $X'$, and for any $x\in X\setminus\fX$ let $x'\in X'$ denote the corresponding point.
We root $X'$ at $x_{\ast}'$.
Then $(X', d', x_{\ast}')$ is a rooted \greg{$\bR$-tree}, which we think of as the `breadth-first spanning tree' of $(X, d, x_\ast)$. 
\greg{Consider the injective map $x\mapsto x'$ defined on $X\setminus\fX$ and let $p:X'\to X$ be the unique continuous extension of its inverse to the whole of $X'$.
Then $p^{-1}(x)=\{x'\}$ for $x\in X\setminus \fX$, and $p^{-1}(x)$ has exactly two elements $x^{(1)},x^{(2)}$ for $x\in \fX$, and both are leaves of $X'$ with $d'(x_{\ast}',x^{(1)})=d'(x_{\ast}',x^{(2)})=d(x_{\ast},x)$.}

If $(X, d, x_\ast, \mu)$ is a rooted measured $\bR$-graph and $\mu$ is non-atomic, then we endow $X'$ with $\mu'$--the natural measure on $X'$ induced by $\mu$\greg{, that is, $\mu'(A)=\mu(p(A))$ for a measurable $A$ in $(X',d')$ and $p$ is as described above.
 This yields a mapping 
\begin{align}\label{eqn:103}
\big(X,d,x_{\ast},\mu\big)\mapsto \big(X',d',x_{\ast}', \{\{x^{(1)},x^{(2)}\} : x\in \fX\},\mu'\big)
\end{align}
from the set of rooted $\bR$-graphs that satisfy \eqref{eqn:1} and are endowed with a non-atomic measure to the set of rooted $\bR$-trees with a non-atomic measure and a finite number of distinguished pairs of leaves that are equidistant to the root. 
This mapping is inverted simply by identifying the distinguished leaves in pairs and taking the metric quotient. We do not give all details since we are not going to use this correspondence explicitly, but it helps in understanding some parts of our constructions. }

The space $X'$ has some nice properties as we describe next.
Note that in going from $X$ to $X'$, distances from the root are preserved: For any $z\in X'$, 
$d'(z, x_{\ast}')=d(p(z), x_{\ast})$, and in particular, 
\begin{align}\label{eqn:2}
\radi(X)=\radi(X').
\end{align}
Further, for any $r>0$,
\begin{align}\label{eqn:4}
\mu\big(B(x_{\ast}, r; X)\big)
=
\mu'\big(B(x_{\ast}', r; X')\big).
\end{align}

For $a<b$, we call a continuous function $e:[a,b]\to[0,\infty)$ with $e(a)=e(b)=0$ a continuous excursion on $[a, b]$.
For a continuous excursion $e$ on $[0,1]$, write $\big(\cT_e, d_e, \rho_e, \mu_e\big)$ for the rooted (at $\rho_e$) measured \greg{$\bR$-tree} encoded by $e$; see Section \ref{sec:constructions} for a precise definition.
Since $X'$ is \greg{an $\bR$-tree} (recall that in our definition, \greg{$\bR$-trees} are always compact), by  \cite[Corollary 1.2]{duquesne2006coding}, there exists a continuous excursion $e_{\bullet}$ on $[0,1]$ such that $(X', d', x_{\ast}')$ is isometric to $\big(\cT_{e_{\bullet}}, d_{e_{\bullet}}, \rho_{e_{\bullet}}\big)$ as rooted metric spaces.
Thus, using \eqref{eqn:2}, we get the following description of the radius of $X$:
\begin{align}\label{eqn:3}
\radi(X)=\|e_{\bullet}\|_{\infty}:=\sup\big\{ e_{\bullet}(t)\ :\ t\in[0,1]\big\}.
\end{align}
If $(X', d', x_{\ast}', \mu')$ is isometric to
$\big(\cT_{e_{\bullet}}, d_{e_{\bullet}}, \rho_{e_{\bullet}},
\mu_{e_{\bullet}}\big)$ as rooted metric measure spaces, then using
\eqref{eqn:4} \chng{and writing $\lambda$ for the Lebesgue measure on $\bR$,} we get
\begin{align}\label{eqn:5}
\mu\big(B(x_{\ast}, r; X)\big)
\chng{
=\mu'\big(B(x_{\ast}', r; X')\big)	
=\mu_{e_\bullet}\big(B(\rho_{e_\bullet},r;\, \cT_{e_\bullet})\big)}=
\lambda\big(\{t\in[0,1]\, :\, e_{\bullet}(t)\leq r\}\big)\, ,
\end{align}
\chng{where the last step follows from the definition of the measured $\bR$-tree coded by $e_{\bullet}$; 
see the beginning of Section~\ref{sec:constructions} for the precise definition.}

\greg{For any $a>0$ and any $f\in C[0, a]$ (not necessarily an excursion), we can define the  occupation measure of $f$ at time $t\in [0,a]$, which we denote by $\fN(f;t,dy)$, by the formula 
\[
\int_0^t \phi(f(u)) du = \int_{y\in \bR}\phi(y)\, \fN(f;t,dy)
\]
for every bounded continuous $\phi:\bR\to\bR$. 
Now suppose that $\fN(f;t,dy)=\eta(f;t,y)dy$ is absolutely continuous with respect to the Lebesgue measure for every $t\in[0, a]$, and that we can choose a version of the density $\eta(f;t,y)$ that is jointly continuous in $t\in [0,a]$ and $y\in \bR$. Then we call $\eta(f;t,y)$ the {\em local time} of $f$ at level $y$ and time $t$, and say that $f$ admits a continuous local time. From this definition it follows that for each fixed $y\in\bR$, the function $\eta(f;\ \cdot\ ,\ y)$ is non-decreasing on $[0,a]$, and the corresponding Stieltjes measures $\eta(f; dt, y)$ satisfy
\begin{align}\label{eqn:6}
\int_0^a \phi\big(t, f(t)\big)dt
=
\int_{y\in\bR}dy\int_{t=0}^a \phi(t, y)\ \eta(f; dt, y)
\end{align}
for every bounded continuous $\phi:[0, a]\times\bR\to\bR$.
}

Now if $(X', d', x_{\ast}', \mu')$ is isometric to $\big(\cT_{e_{\bullet}}, d_{e_{\bullet}}, \rho_{e_{\bullet}}, \mu_{e_{\bullet}}\big)$ for a continuous excursion $e_{\bullet}$ on $[0, 1]$ and $e_{\bullet}$ admits a continuous local time, then the right side of \eqref{eqn:5} is differentiable w.r.t. $r$, and further,
\begin{align}\label{eqn:7}
\frac{d}{d r}\mu\big(B(x_{\ast}, r; X)\big)
=
\eta(e_{\bullet}; 1, r)\ \ \text{ for all }\ \ r>0\, .
\end{align}
The left side of \eqref{eqn:7} gives the distance profile around the root $x_\ast$.

Thus, if $X$ is a random rooted measured $\bR$-graph and we can identify the random excursion $e_{\bullet}$ that encodes $X'$, \greg{and $e_\bullet$ admits a continuous local time}, then we can express the radius of $X$ and the distance profile around the root of $X$ in terms of the supremum of $e_{\bullet}$ and the local time of $e_{\bullet}$.
For the random $\bR$-graphs considered in this paper, the root $x_{\ast}$ will be a $\mu$-distributed point in $(X, d, \mu)$. 
In this case the two-point function of $(X, d, \mu)$ will have the same law as $e_{\bullet}(U)$, where $U\sim\mathrm{Uniform}[0, 1]$ is independent of $e_{\bullet}$.

In this paper we will consider the following two random $\bR$-graphs:

\begin{inparaenumi}
\noindent\item 
the scaling limit $\cH_{(s)}$ of uniform connected rooted (labeled) graphs with surplus $s$, and
\noindent\item
the scaling limit $\crumg$ of random unicellular maps of genus $g$ \greg{considered in \cite{lab-oa-gc-ef-cg}}.
\end{inparaenumi}

\noindent
The precise definitions of these spaces will be given in Section \ref{sec:constructions}.
We will identify the distributions of their breadth-first spanning trees $\cH_{(s)}'$ and $\crumg'$, and we will describe how to glue random points on these \greg{$\bR$-trees} to recover the spaces  $\cH_{(s)}$ and $\crumg$ in distribution.
%As described above, these constructions will in particular enable us to express the radii, distance profiles, and two-point functions in $\cH_{(s)}$ and $\crumg$ in terms of functionals of suitably tilted Brownian excursions.

The spaces $\cH_{(s)}$, $s\geq 0$, are the building blocks for the scaling limit of the critical \erdos random graph identified in \cite{BBG-12}; 
see \cite[Construction 3.10 and Theorem 3.12]{addarioberry-sen} for a precise statement.
The \erdos scaling limit is a universal object in the sense that it arises as the scaling limit of a wide array of standard models of critical random discrete structures exhibiting mean-field behavior. 
Examples of such models include random regular graphs under critical percolation or more generally critical random graphs with given degree sequence (under finite third moment assumption on the degrees), various models of inhomogeneous random graphs (under appropriate assumptions), bounded-size rules, and the vacant set left by random walks on random regular graphs; see \cite{SBSSXW14, bhamidi-broutin-sen-wang, bhamidi-sen}.
Further, existing literature suggests that the components of the high-dimensional discrete torus \cite{hofstad-sapozhnikov,heydenreich-hofstad-1,heydenreich-hofstad-2} and the hypercube \cite{hofstad-nachmias} under critical percolation, and the critical quantum random graph model \cite{dembo-levit-vadlamani} also share the \erdos scaling limit.
The breadth-first construction of $\cH_{(s)}$ (given in Theorem \ref{thm:bf=df-Hs} below) in particular yields the same for this universal scaling limit.

The construction of the \erdos scaling limit given in \cite{BBG-12} can be thought of as a `depth-first construction'; see Construction \ref{constr:H-s} below.
An alternate construction that can be seen as a `core decomposition' was given in \cite{BBG-limit-prop-11}.
The advantage of the breadth-first perspective, as explained earlier, is that it directly identifies the radius, the two-point function, and the distance profile in terms functionals of suitably tilted Brownian excursions, which the two constructions above do not.
The problem of identifying the breadth-first construction of the \erdos scaling limit was asked in \cite[Section 6]{BBG-12}.
This was the main motivation behind this work.
Quoting the authors of \cite{BBG-12}:
\begin{quote}
	The rescaled breadth-first walk $\ldots$ converges to the same limit as the rescaled {\it height profile} (i.e. the number of vertices at each height) of a ``breadth-first tree", which contains less
	information and, in particular, does not code the structure of that tree. As a result, it
	seems that it would be much harder to derive a metric space construction of a limiting
	component using the breadth-first viewpoint. 
\end{quote}
The breadth-first approach indeed requires more than just the convergence of the rescaled breadth-first walk.
This will be discussed in more detail in Section \ref{sec:discussion}.

\medskip

\noindent{\bf Comments on the use of maps in our approach:} 
\greg{Our proofs will make a thorough use of maps, %i.e.\ embedded graphs in surfaces, 
not only in the study of $\crumg$, but also in that of the spaces $\cH_{(s)}$. This might look a bit surprising at first sight. Indeed, the spaces $\cH_{(s)}$ were initially introduced as scaling limits of random labeled graphs. The most natural approach would be to perform a breadth-first search on the same family of graphs, cutting cycles as they appear in this exploration in order to obtain a random labeled tree with a certain law that one can compare to that of a uniform rooted labeled tree $H_{n, 0}$ of a given size $n$. Running the construction backwards, to obtain the random labeled graph starting from this random tree, one has to choose pairs of vertices in the tree at (roughly) the same height and connect them by an edge. A difficulty then appears when trying to understand the scaling limit of this construction: the operation of choosing random pairs of vertices at the same height in a continuum tree coded by some excursion $e$ amounts to sampling them according to measures associated with the local time of $e$, and proving the convergence then requires a very good understanding of the discrete versions of the local time associated with either the height function or the contour function of the tree $H_{n, 0}$ as $n\to\infty$. Unfortunately, such results are not available in the literature on random trees.
%(the relevant case in our study would be Galton-Watson trees with Poisson(1) offspring distribution, conditioned on the total population). 
}

\greg{Instead, we take an indirect approach by considering a different model: uniform maps with a fixed surplus.
Using Construction \ref{constr:H-s} given below, we will show that the scaling limit of this model is also $\cH_{(s)}$.
Now, by considering a breadth-first exploration of this model, we circumvent the issue with local times mentioned above by connecting the original problem to a problem about suitably tilted plane trees (as opposed to labeled trees).
The contour function of a uniform plane tree is a simple random walk excursion, and this will allow us to apply existing results about local time fields of random walks in our proofs.
Owing to this reason, the proof of the breadth-first construction of $\cH_{(s)}$ proceeds via a study of such maps.}
%This will also enable us to deduce the results for $\crumg$ easily.
\medskip

\noindent{\bf Organization of this paper:} The rest of this paper is organized as follows:
In Section \ref{sec:constructions}, we describe the constructions of the random metric spaces involved.
In Section \ref{sec:main-result}, we state our main results.
In Section \ref{sec:notation}, we explain the notation and conventions used in this paper; in particular, the notation related to maps will be described.
The proofs of our main results will be given in various subsections in Section \ref{sec:proofs}.
Section \ref{sec:discussion} contains some related questions and further discussions.
The proof of a key result used in Section \ref{sec:proofs} will be outlined in Appendix \ref{sec:appendix}.

\medskip

\noindent{\bf Convention about metric spaces:}
We fix a convention here that we will follow throughout the paper.
For any metric measure space $\mvX = (X, d, \mu)$ and $\alpha>0$, $\alpha\cdot\mvX$ or simply $\alpha\mvX$ will denote
the metric measure space $(X, \alpha d, \mu)$, i.e, the space where the metric has been multiplied
by $\alpha$ and the measure $\mu$ has remained unchanged.
We can similarly define $\alpha\mvX$ for a rooted metric measure space $\mvX$ by leaving the root unchanged.

When dealing with convergence of rooted metric measure spaces, we will work with the topology induced by the `pointed Gromov-Hausdorff-Prokhorov (GHP) distance,' which we denote by $d_{\GHP}^1(\cdot,\cdot)$.
We refer the reader to \cite[Section 2.1]{AddBroGolMie13} for the relevant definitions.

\section{Constructions of the spaces}\label{sec:constructions}
Let $e$ be a continuous excursion on $[0, 1]$.
Let $d_e$ be the pseudo-metric on $[0,1]$ given by
\begin{equation}\label{eqn:8}
d_e(s,t):= e(s) + e(t) - 2 \inf_{u \in [s,t]}e(u), \; \mbox{ for } s,t  \in [0,1].
\end{equation}
Define the equivalence relation $s \sim t \Leftrightarrow d_e(s,t) = 0$. 
Let $[0,1]/\sim$ denote the corresponding quotient space and consider the space $\cT_e:= [0,1]/\sim$
endowed with the quotient metric on the equivalence classes induced by $d_e$. 
We abuse notation and write $d_e$ for the quotient metric on $\cT_e$ as well.
Then $(\cT_e, d_e)$ is \greg{an $\bR$-tree} \cite{legall-survey,evans-book}.
Let $q_e:[0,1] \to \cT_e$ be the canonical projection and write $\mu_e$ for the push-forward of the Lebesgue measure on $[0,1]$ onto $\cT_e$ via $q_e$. 
Further, we let $\cT_e$ be rooted at $\rho_e := q_e(0)$.
Then $(\cT_e, d_e, \rho_e, \mu_e)$ is a rooted measured \greg{$\bR$-tree}.
Note that by construction, for any $x\in \cT_e$, the function $e$ is constant on $q_e^{-1}(\{x\})$.
%Thus for each $x\in [0,l]$, we write $\hght(x) = h(q_h^{-1}(x))$ for the height of this vertex.
Note also that the height of $x\in\cT_e$ defined as $\hght(x; \cT_e):=d_e(\rho_e, x)$ satisfies 
$\hght(x; \cT_e)=e(u)$ for any $u\in q_e^{-1}(\{x\})$.
We define the set of leaves of $\cT_e$ to be
\[
\cL(\cT_e)
:=
\big\{x\in\cT_e\, :\, \degr(x; \cT_e)=1  \}.
\]

We will write $\big(\ve(t)\, ,\, t\in[0, 1]\big)$ for a standard Brownian excursion.
The \greg{$\bR$-tree} $2\cdot\cT_{\ve}$ is called the Brownian
continuum random tree
\chng{(recall our convention about metric spaces explained at the end of Section~\ref{sec:intro}).}
It is well-known \cite{aldous-crt-1, aldous-crt-3} that the measure $\mu_{\ve}$ (also called the mass measure) on $\cT_{\ve}$ is non-atomic and concentrated on $\cL(\cT_{\ve})$ almost surely.

\subsection{Constructions of the scaling limit of uniform connected graphs with fixed surplus}

We will now define the random spaces $\cH_{(s)}$, $s\geq 0$, that were introduced in Section~\ref{sec:intro}.
For a finite connected graph $G=(V, E)$, let $\spls(G):=|E|-|V|+1$ denote the number of surplus \greg{(also called excess)} edges in $G$.
For $n\geq 1$ and $s\geq 0$ define
\begin{align}\label{eqn:64}
\bH_{n, s}=\big\{G\, :\, G\text{ rooted, connected, simple, labeled graph on }[n]\text{ with }\spls(G)=s \big\}\, ,
\end{align}
where $[n]:=\{1,\ldots, n\}$.
Let $H_{n, s}$ be uniformly distributed over $\bH_{n, s}$.
View $H_{n, s}$ as  a rooted metric measure space by endowing it with the graph distance and the uniform probability measure on the vertices.
Then there exists a random compact, rooted metric measure space $\cH_{(s)}$ such that
\begin{align}\label{eqn:9}
n^{-1/2}\cdot H_{n, s}\weakc\cH_{(s)}
\end{align}
w.r.t. the pointed GHP topology.
The space $\cH_{(0)}$ is simply $2\cdot \cT_{\ve}$--the Brownian continuum random tree \cite{aldous-crt-1, aldous-crt-3}.
For $s\geq 1$, \eqref{eqn:9} can be proved by arguments similar to the ones used in  \cite{BBG-12, bhamidi-sen};
a brief sketch of the proof is given in \cite[Section A.1]{addarioberry-sen}.

The spaces $\cH_{(s)}$, $s\geq 0$, are central to describing the scaling limits of many other random discrete structures.
As mentioned before, the critical \erdos scaling limit identified in \cite{BBG-12} can be expressed in terms of the spaces $\cH_{(s)}$, $s\geq 0$.
The spaces $\cH_{(s)}$ also arise as the scaling limit of

\noindent (i) uniform connected graphs with a given degree sequence under some assumptions on the degree sequence \cite[Theorem 2.4]{bhamidi-sen}, and 

\noindent (ii) uniform rooted maps with fixed surplus $s$ (without any restriction on the genus); see Remark \ref{rem:uniform-map-scaling-limit} below. 

\noindent The scaling limit of the minimal spanning tree of the complete graph identified in \cite{AddBroGolMie13} can be expressed as the limit, as $s\to\infty$, of the space obtained by applying a `cycle-breaking' procedure on the space $(12s)^{1/6}\cdot\cH_{(s)}$ \cite[Theorem 4.8]{addarioberry-sen}.
We now describe one construction of $\cH_{(s)}$.

\begin{constr}[Depth-first construction of $\cH_{(s)}$]\label{constr:H-s}
	Fix an integer $s\geq 0$.
	\begin{enumeratea}
		\item Sample $\esdf$ with law given by
		\begin{align}\label{eqn:13}
		\bE\big[\phi(\esdf)\big]=\frac{\bE\big[\phi(\ve)\big(\int_0^1 \ve(t)dt\big)^s\big]}{\bE\big[\big(\int_0^1 \ve(t)dt\big)^s\big]}\,
		\end{align}
		for every bounded continuous $\phi:C[0, 1]\to\bR$.
		\item Conditional on $\esdf$, sample i.i.d. points $u_1,\ldots, u_s$ in $[0, 1]$ with density
		\[\left(\frac{\esdf(u)}{\int_0^1 \ve_{(s)}^{\df}(t)dt}\right)du,\ \ u\in[0, 1].\]
		\item
		Conditional on $\esdf$ and $u_1,\ldots, u_s$, sample independent points $z_1,\ldots,z_s$, where $z_i$ is uniformly distributed in $[0, \esdf(u_i)]$.
		For $1\leq i\leq s$, let
		\[
		v_i:=\inf\big\{t\geq u_i\ :\ \esdf(t)=z_i\big\}.
		\]
		\item Set $\cH_{(s)}$ to be the quotient space $2\cdot\big(\cT_{\esdf}\; /\sim\big)$, where $\sim$ is the smallest equivalence relation under which $q_{\esdf}(u_i)\sim q_{\esdf}(v_i)$ for $1\leq i\leq s$. \chng{Here, $\cH_{(s)}$ is endowed with the root and the measure inherited from $\cT_{\esdf}$ by the quotient map}.
		That is, $\cH_{(s)}$ is the rooted measured $\bR$-graph obtained by identifying the points $q_{\esdf}(u_i)$ and $q_{\esdf}(v_i)$ on $\cT_{\esdf}$ and then multiplying the distances by $2$.
	\end{enumeratea}
\end{constr}
The above construction of the space $\cH_{(s)}$ is essentially
contained in the arguments given in \cite{BBG-12}. 
The reason for using the notation $\esdf$ is that Construction
\ref{constr:H-s} can be thought of as the continuum analogue of the
depth-first construction of a finite connected graph. 
The \greg{$\bR$-tree} $2\cdot\cT_{\esdf}$ plays the role of the
depth-first spanning tree of $\cH_{(s)}$.

\chng{
One can heuristically argue that the expression for the tilt in \eqref{eqn:13} has to be proportional to $\big(\int_0^1 \ve(t)dt\big)^s$ as follows:
The space $\cH_{(s)}$ can be seen as a ``uniform rooted measured $\bR$-graph having $s$ cycles" (see \eqref{eqn:9}) just as the Brownian continuum random tree $2\cdot\cT_{\ve}$ is, in a sense, a ``uniform random rooted measured $\bR$-tree."
Construction~\ref{constr:H-s} gives a construction of $\cH_{(s)}$ by identifying $s$ i.i.d. pairs of points on $2\cdot\cT_{\esdf}$, where each pair is uniformly distributed subject to the constraint that the pair consists of a point and an ancestor of that point.
Accordingly, the law of $2\cdot\cT_{\esdf}$ is the one obtained by tilting the law of $2\cdot\cT_{\ve}$ by the total ``weight'' of possible chosen pairs, which equals
\[
\Big(\int_{x\in \cT_{\ve}}2\cdot\hght(x; \cT_{\ve})\, \mu_{\ve}(dx)\Big)^s
=
\Big(\int_0^1 2\ve(t) dt\Big)^s\, .
\]
Alternately, one can arrive at the expression for the tilt by looking at the analogous problem in the discrete setting; see the proof of Proposition~\ref{prop:df-scaling-limit}.
}

\greg{We will now define another space $\hsbf$.
Recall, as is well-known, that $\ve$ admits a.s. a continuous local time $\eta(\ve;\cdot,\cdot)$ as defined  around \eqref{eqn:6}.
The same is true of any random process with a law that is absolutely continuous with respect to that of $\ve$, which justifies the following construction.}

\begin{constr}[The space $\hsbf$]\label{constr:H-s-bf}
	Fix an integer $s\geq 0$. 
	\begin{enumeratea}
		\item
		 Sample $\esbf$ with law given by
		\begin{align}\label{eqn:14}
		\bE\big[\phi(\esbf)\big]
		=
		\frac{\bE\Big[\phi(\ve)\cdot\big(\int_0^{\infty}\eta(\ve; 1, y)^2 dy\big)^s\Big]}{\bE\Big[\big(\int_0^{\infty}\eta(\ve; 1, y)^2 dy\big)^s\Big]}
		\end{align}
		for every bounded continuous $\phi:C[0, 1]\to\bR$.
		Write 
		\begin{align}\label{eqn:38}
		\nsbf(\cdot, \cdot)=\eta\big(\esbf\,;\, \cdot, \cdot\big).
		\end{align}
		\item
		Conditional on $\esbf$, sample i.i.d. points $H_1,\ldots, H_s$ in $[0,\infty)$ with density
		\begin{align}\label{eqn:42}
		\left(\frac{\nsbf(1, h)^2}{\int_0^{\infty}\nsbf(1, y)^2 dy}\right)dh.
		\end{align}
		\item
		Conditional on $\esbf$ and $H_1,\ldots, H_s$, sample $2s$ independent points $u_1, v_1,\ldots,u_s, v_s$ in $[0, 1]$, where $u_i$ and $v_i$ are distributed according to \greg{the} law
		\[\frac{\nsbf\big(dt, H_i\big)}{\nsbf\big(1, H_i\big)}\]
		for $1\leq i\leq s$.
		\item
		Set $\hsbf$ to be the quotient space
                $2\cdot\big(\cT_{\esbf}\; /\sim\big)$, where $\sim$ is the smallest equivalence relation under which $q_{\esbf}(u_i)\sim q_{\esbf}(v_i)$ for $1\leq i\leq s$.
                \chng{Here, $\hsbf$ is endowed with the root and the measure inherited from $\cT_{\esbf}$ by the quotient map.} 
		That is, $\hsbf$ is the rooted measured $\bR$-graph
                obtained by identifying the points $q_{\esbf}(u_i)$
                and $q_{\esbf}(v_i)$ on $\cT_{\esbf}$ and then
                multiplying the distances by $2$. 
	\end{enumeratea}
\end{constr}

Let us describe the above construction in words.
We first sample a tilted excursion $\esbf$.
Then we sample $s$ i.i.d. `heights' $H_1,\ldots, H_s$ according to
density \eqref{eqn:42}. 
For $1\leq i\leq s$, we sample two points in the tree $\cT_{\esbf}$ independently according to the normalized local time measure at height $H_i$, and then we glue these points.
Finally, we multiply the metric by $2$ in the resulting
space.

Almost surely, for every $y>0$ the measure $\eta\big(\ve; dt, y\big)$ is concentrated on $q_{\ve}^{-1}\big(\cL(\cT_{\ve})\big)$.
By absolute continuity, the same is true if we replace $\ve$ by $\esbf$.
Thus $q_{\esbf}(u_i)$ and $q_{\esbf}(v_i)$ are leaves in $\cT_{\esbf}$, $1\leq i\leq s$.
Further, almost surely $\esbf(u_i)=\esbf(v_i)=H_i$, $1\leq i\leq s$, and consequently, $q_{\esbf}(u_i)$ and $q_{\esbf}(v_i)$ are equidistant to the root of $\cT_{\esbf}$.
\greg{So we see that $\hsbf$ is obtained from $2\cdot \cT_{\esbf}$ and the points $\{q_{\esbf}(u_i),q_{\esbf}(v_i)\}$, $1\leq i\leq s$, by applying the inverse of the operation described in \eqref{eqn:103}. }
%So if $w_i$ is the point in $\hsbf$ resulting from the identification of $q_{\esbf}(u_i)$ and $q_{\esbf}(v_i)$, then $\deg\big(w_i;\, \hsbf\big)=2$, $1\leq i\leq s$.
%Since $\hght\big(u_i; \cT_{\esbf}\big)=\hght\big(v_i; \cT_{\esbf}\big)=H_i$, $1\leq i\leq s$, the cut locus of the root in $\hsbf$ is $\big\{w_1,\ldots, w_s\big\}$ and $\hsbf$ satisfies \eqref{eqn:1}.
Hence the breadth-first spanning tree of $\hsbf$ is $(\hsbf)'=2\cdot\cT_{\esbf}$.

Theorem \ref{thm:bf=df-Hs} given below states that $\cH_{(s)}\equald\hsbf$.
Thus, Construction~\ref{constr:H-s-bf} gives the breadth-first construction of the space $\cH_{(s)}$.
\chng{
	Similar to Construction~\ref{constr:H-s}, one can arrive at the expression for the tilt in \eqref{eqn:14} using heuristic arguments.
	As mentioned before, the space $\cH_{(s)}$ can be seen as a uniform rooted measured $\bR$-graph having $s$ cycles, and the aim in Construction~\ref{constr:H-s-bf} is to construct a new space $\hsbf$ in a `breadth-first fashion' so that we have $\cH_{(s)}\equald\hsbf$ .
	Now, in Construction~\ref{constr:H-s-bf} we have constructed $\hsbf$ by identifying $s$ i.i.d. pairs of points on $2\cdot\cT_{\esbf}$, where each pair is uniformly distributed subject to the constraint that the points in the pair are at the same height.
	Note that the push-forward of the Lebesgue measure on $[0, 1]$ onto $[0, 1]\times\bR$ under the mapping $t\mapsto (t, \ve(t))$ is $\eta(\ve; dt, h)dh$, and accordingly, the law of $2\cdot\cT_{\esbf}$ is the one obtained by a tilting of the law of $2\cdot\cT_{\ve}$ by the total ``weight'' of possible chosen pairs, which equals
	\[
	\Big(\int_{h=0}^{\infty}dh\big(\int_{u=0}^1\int_{v=0}^1 \eta(\ve; du, h)\eta(\ve; dv, h)\big)\Big)^s
	=
	\big(\int_0^{\infty}\eta(\ve; 1, h)^2 dh\big)^s\, .
	\]
	Alternately, one can arrive at the expression for the tilt by looking at the analogous problem in the discrete setting; see the proof of Proposition~\ref{prop:bf-scaling-limit}.
}

\begin{rem}\label{rem:not-same-distribution}
It follows from Jeulin's local time identity \cite{AMP, jeulin} that
\begin{align}\label{eqn:74}
\int_{0}^{\infty}\eta\big(\ve; 1, y\big)^2 dy\,
\equald\,
2\int_0^1\ve(t)dt\, ;
\end{align}
see the argument given around \eqref{eqn:101}.
However $\esdf$ and $\esbf$, as defined in \eqref{eqn:13} and \eqref{eqn:14} respectively, do not have the same distribution.
In fact, from the discussion above and \eqref{eqn:3}, $\radi\big(\cH_{(s)}\big)\equald 2\cdot\|\esbf\|_{\infty}$.
However, from Construction \ref{constr:H-s}, it is clear that $\radi\big(\cH_{(s)}\big)$ is stochastically dominated by $2\cdot\|\esdf\|_{\infty}$.
\end{rem}

\subsection{Constructions of the continuum random unicellular maps}\label{sec:def-S-g}

Fix $g\geq 1$, and let $\bum$ be the set of rooted unicellular maps of genus $g$ having $n+1$ vertices.
Let $UM_{n, g}$ be uniformly distributed over $\bum$.
Denote its root edge by $e_{\ast}$.
As before, we endow $UM_{n, g}$ with the graph distance and the uniform probability measure on the vertices, and think of it as a rooted metric measure space with the root being $e_{\ast}(0)$--the origin of $e_{\ast}$. 
(The notation related to maps will be discussed in Section \ref{sec:notation}.)
Then there exists a random compact rooted metric measure space $\crumg$ such that
\begin{align}\label{eqn:10}
n^{-1/2}\cdot UM_{n, g}\weakc\crumg
\end{align}
w.r.t. the pointed GHP topology.
As mentioned in \cite[Page 940]{lab-oa-gc-ef-cg}, it seems that a proof of the convergence \eqref{eqn:10} is not written down in the literature.
However, the result can be deduced by following the arguments used in \cite{BBG-limit-prop-11} or in \cite{BBG-12}.
A construction of $\crumg$ that can be viewed as a `core decomposition' is given in \cite{lab-oa-gc-ef-cg}\footnote{In fact, the authors consider a more general model in \cite{lab-oa-gc-ef-cg}.}.
We will next describe the breadth-first construction of $\crumg$.
To do so, we first need to set up some notation.

Let $\bS_{(g)}$ be the set of permutations on $[4g]=\{1, 2,\ldots, 4g\}$ that satisfy the following:
$\sigma\in \bS_{(g)}$ iff 
\begin{enumeratei}
\item 
the cycle decomposition of $\sigma$ has $2g$ many transpositions, and
\item 
the permutation $\varrho\sigma$ on $[4g]$ has only one cycle of length $4g$, where $\varrho=\big(1, 2, \ldots, 4g\big)$.
\end{enumeratei}
For example, $\bS_{(1)}=\big\{(1, 3)(2, 4)\big\}$, and $(1, 3)(2, 4)(5, 7)(6, 8)$ and $(1, 7)(2, 5)(3, 8)(4, 6)\in \bS_{(2)}$.

Suppose $\sigma=\prod_{j=1}^{2g}(\ell_{2j-1}, \ell_{2j})\in \bS_{(g)}$ with 
$\ell_1<\ell_3<\ldots<\ell_{4g-1}$, and $\ell_{2j-1}<\ell_{2j}$ for each $j\in [2g]$.
For a continuous excursion $e$ on $[0, 1]$ that admits a continuous local time (as discussed around \eqref{eqn:6}), and $\mvy:=(y_1,\ldots, y_{2g})\in (0, \infty)^{2g}$, 
define a measure $\nu_{e, \sigma, \mvy}$ on $[0, 1]^{4g}$ as follows:
\begin{align}\label{eqn:12}
\nu_{e, \sigma, \mvy}\big(dt_1,\ldots, dt_{4g}\big)
:=
\ind\{t_1<\ldots<t_{4g}\}\cdot
\prod_{j=1}^{2g}\Big(\eta\big(e; dt_{\ell_{2j-1}}, y_j\big)\eta\big(e; dt_{\ell_{2j}}, y_j\big)\Big)\, .
\end{align}
The push-forward of the probability measure $\nu_{e, \sigma, \mvy}(\cdot)/\nu_{e, \sigma, \mvy}([0, 1]^{4g})$ onto $(\cT_e)^{4g}$ under the $4g$-fold product of the quotient map $q_e$ can be viewed as the `uniform measure' on the set 
$\big\{\big(q_e(t_1),\ldots, q_e(t_{4g})\big)\, :\, 
0<t_1<\ldots< t_{4g}<1\text{ and }
e\big(t_{\ell_{2j-1}}\big)=e\big(t_{\ell_{2j}}\big)=y_j\text{ for all }j\in[2g]\big\}.$

\begin{constr}[The space $\crumgbf$]\label{constr:crum-g-bf}
	Fix an integer $g\geq 1$.
	\begin{enumeratea}
		\item
		Sample $\egum$ with law given by
		\begin{align}\label{eqn:78}
		\bE\big[\phi(\egum)\big]
		=
		24^g\cdot g!\cdot
		\bE\Big[\phi(\ve)\cdot
	    \sum_{\sigma\in \bS_{(g)}}
	    \Big(
	    \int_{(0,\infty)^{2g}}
	    \nu_{\ve, \sigma, \mvy}\big([0, 1]^{4g}\big)\ dy_1\ldots dy_{2g}
	    \Big)
	    \Big]
		\end{align}
		for every bounded continuous $\phi:C[0, 1]\to\bR$.
		(That this indeed gives a valid probability distribution will be shown in \eqref{eqn:98}.)
		\item
		Conditional on $\egum$, sample a permutation $\Theta$ with distribution
		\[
		\pr\big(\Theta=\theta\big)=
		\frac{\int_{(0,\infty)^{2g}}
			\nu_{\egum, \theta, \mvy}\big([0, 1]^{4g}\big)\ dy_1\ldots dy_{2g}
			}{\sum_{\sigma\in \bS_{(g)}}\Big(
			\int_{(0,\infty)^{2g}}
			\nu_{\egum, \sigma, \mvy}\big([0, 1]^{4g}\big)\ dy_1\ldots dy_{2g}
			\Big)}\, ,\ \ \ 
		\theta\in \bS_{(g)}\, .
		\]
		\item
		Conditional on $\egum$ and $\Theta$, sample $\mvH=(H_1,\ldots, H_{2g})$ with density
		\[
		\frac{\nu_{\egum, \Theta, \mvh}\big([0, 1]^{4g}\big)\, dh_1\ldots dh_{2g}}{\int_{(0,\infty)^{2g}}
			\nu_{\egum, \Theta, \mvy}\big([0, 1]^{4g}\big)\ dy_1\ldots dy_{2g}},\ \ \ 
		\mvh=(h_1,\ldots, h_{2g})\in(0,\infty)^{2g}\, .
		\]
		\item 
		Conditional on $\egum$, $\Theta$, and $\mvH$, sample $(u_1,\ldots, u_{4g})$ according to \greg{the} law
		\[
		\frac{\nu_{\egum, \Theta, \mvH}\big(\cdot\big)}{\nu_{\egum, \Theta, \mvH}\big([0, 1]^{4g}\big)}\, .
		\]
		\item
		Set $\, \crumgbf$ to be the quotient space $\cT_{\egum}\; /\sim$, where $\sim$ is the \greg{smallest} equivalence relation under which $q_{\egum}(u_i)\sim q_{\egum}(u_{\sigma(i)})$ for $ i\in [4g]$.
		Thus, $\crumgbf$ is the rooted measured $\bR$-graph obtained by identifying the points $q_{\egum}(u_i)$ and $q_{\egum}(u_{\sigma(i)})$ on $\cT_{\egum}$ for all $i$.
	\end{enumeratea}
\end{constr}

The construction of $\crumgbf$ has more information than just the metric measure space structure--$\crumgbf$ can in fact be viewed as a continuum map; see \cite{lab-oa-gc-ef-cg} for a more detailed discussion on continuum random maps.
The defining conditions for $\bS_{(g)}$ ensure that the cycles resulting from the identifications $q_{\egum}(u_i)\sim q_{\egum}(u_{\sigma(i)})$ are appropriately `entangled' so that the resulting map is unicellular; see \cite{chapuy} for a detailed account of structure of unicellular maps.

By arguments similar to the ones given below Construction \ref{constr:H-s-bf}, 
$q_{\egum}(u_i)$ is a leaf in $\cT_{\egum}$ for each $i\in [4g]$.
Further, 
$\hght\big(q_{\egum}(u_i);\, \cT_{\egum}\big)
=
\hght\big(q_{\egum}(u_{\sigma(i)});\, \cT_{\egum}\big)$ 
for all $i$.
It thus follows that $\big(\crumgbf\big)'=\cT_{\egum}$.
Theorem \ref{thm:bf=df-crum} given below states that $\crumg\equald\crumgbf$.
Thus, Construction \ref{constr:crum-g-bf} gives the breadth-first construction of $\crumg$.

\begin{rem}
From Construction \ref{constr:H-s-bf}, the space $\hsbf$ can also be viewed as a continuum map.
Then it can be shown that the space $\crumgbf$ as given in Construction \ref{constr:crum-g-bf} has the same law as `$\htwogbf$ conditioned to be unicellular.'
More precisely, let $u_i, v_i$, $1\leq i\leq 2g$, be as in Construction \ref{constr:H-s-bf} (corresponding to $s=2g$), and order them as $t_1<\ldots<t_{4g}$.
Define a permutation $\psi$ on $[4g]$ by letting $\psi(i)=j$ iff $\{t_i, t_j\}=\{u_k, v_k\}$ for some $k\in [2g]$.
Then it can be checked that conditional on $\big\{\psi\in\bS_{(g)}\big\}$, $\htwogbf$ viewed as a continuum map has the same law as $\crumgbf$.
\end{rem}

\section{Main results}\label{sec:main-result}
The following theorem gives the breadth-first construction of the space $\cH_{(s)}$, and consequently, of the critical \erdos scaling limit.

\begin{thm}\label{thm:bf=df-Hs}
Fix $s\geq 0$, and let $\cH_{(s)}$ and $\hsbf$ be as in \eqref{eqn:9} and Construction \ref{constr:H-s-bf} respectively.
Then $\cH_{(s)}\equald\hsbf$.
\end{thm}
As explained in Section \ref{sec:intro}, the following results are immediate from Theorem \ref{thm:bf=df-Hs}.

\begin{cor}\label{cor:bf=df-Hs}
\noindent {\upshape (i) Radius:} We have, $\radi(\cH_{(s)})\equald 2\cdot\|\esbf\|_{\infty}$.
\vskip3pt

\noindent {\upshape (ii) Two-point function:} Denote the metric in $\cH_{(s)}$ by $d(\cdot, \cdot)$  and the measure on $\cH_{(s)}$ by $\mu$.
Let $x_1$ and $x_2$ be two i.i.d. $\mu$-distributed points in $\cH_{(s)}$.
Then
\[
d(x_1, x_2)\equald 2\cdot\esbf(U)\, ,
\]
where $U\sim\mathrm{Uniform}[0, 1]$ and is independent of $\esbf$.
\vskip3pt

\noindent {\upshape (iii) Distance profile:} Denoting the root of $\cH_{(s)}$ by $x_{\ast}$, we have,
\[
\Big(\frac{d}{dr}\mu\big(B\big(x_{\ast}, r; \cH_{(s)}\big)\big)\, ,\ r>0\Big)
\equald
\Big(\eta\big(2\esbf; 1, r\big)\, ,\ r>0\Big)
=
\frac{1}{2}\Big(\eta\big(\esbf; 1, r/2\big)\, ,\ r>0\Big)\, .
\]
\end{cor}

Recall the definition of $\esdf$ from \eqref{eqn:13}.
Using the breadth-first view point we get another representation of $\radi(\cH_{(s)})$ in terms of $\esdf$ which we state in the next corollary.
\begin{cor}\label{cor:radius-Hs}
For any $s\geq 1$,
\[
\radi\big(\cH_{(s)}\big)
\equald
\int_0^1\frac{1}{\esdf(t)}dt\, .
\]
\end{cor}
The corresponding result for $s=0$ is well-known; it says that the height of the Brownian continuum random tree $2\cT_{\ve}$ has the same distribution as $\int_0^1 dt/\ve(t)$.
Corollary \ref{cor:radius-Hs} gives the analogue of this result for graphs.

Note that \eqref{eqn:9} together with Corollary \ref{cor:bf=df-Hs} (i) and (ii) imply that
\begin{align}\label{eqn:15}
n^{-1/2}\cdot\radi\big(H_{n, s}\big)\weakc 2\cdot\|\esbf\|_{\infty}\ \ \
\text{ and }\ \ \
n^{-1/2}\cdot d_n(v_1, v_2)\weakc 2\cdot\esbf(U)\, ,
\end{align}
where $d_n(\cdot, \cdot)$ denotes the graph distance in $H_{n, s}$, and $v_1, v_2$ are independent and uniformly distributed in $[n]$.
However, \eqref{eqn:9} and Corollary \ref{cor:bf=df-Hs} (iii) do not imply the convergence of the (properly rescaled) distance profile in $H_{n, s}$.
The following theorem says that this convergence holds as well.

\begin{thm}\label{thm:H-n-s-width-convergence}
	Fix $s\geq 1$.
	Let $Z_{n, s}(\ell)$ denote the number of vertices in $H_{n, s}$ at distance $\ell$ from the root, $\ell=0,1,\ldots$.
	Let $\bar Z_{n, s}(r)=n^{-1/2}\cdot Z_{n,s}\big(\lfloor r\sqrt{n}\rfloor\big)$, $r\geq 0$.
	Then as $n\to\infty$,
	\[
	\Big(\bar Z_{n, s}(r),\ r\geq 0\Big)
	\weakc
	\Big(\eta\big(2\esbf; 1, r\big)\, ,\ r>0\Big)
	=
	\Big(\frac{1}{2}\eta\big(\esbf; 1, r/2\big),\ r\geq 0\Big)
	\]
	w.r.t. Skorohod $J_1$ topology on $\bD([0,\infty)\ :\ \bR)$.
\end{thm}

The analogue of Theorem \ref{thm:H-n-s-width-convergence} for $s=0$ deals with convergence of the height profiles of uniform rooted labeled trees.
This result, in a more general form, was put forward as a conjecture in \cite[Conjecture 4]{aldous-crt-2} and proved in \cite[Theorem 1.1]{drmota-gittenberger} (see Theorem \ref{thm:total-local-time} below).
Theorem \ref{thm:H-n-s-width-convergence} gives a generalization of this result for graphs.

The next theorem gives the breadth-first construction of the space $\crumg$.

\begin{thm}\label{thm:bf=df-crum}
Fix $g\geq 1$, and let $\crumg$ and $\crumgbf$ be as in \eqref{eqn:10} and Construction \ref{constr:crum-g-bf} respectively.
Then $\crumg\equald\crumgbf$.
\end{thm}
The following result is the analogue of Corollary \ref{cor:bf=df-Hs}.
\begin{cor}\label{cor:bf=df-crum}
The conclusions of Corollary \ref{cor:bf=df-Hs} continue to hold if we replace $\cH_{(s)}$ by $\crumg$ and $2\cdot\esbf$ by $\egum$.
\end{cor}

\section{Notation and conventions}\label{sec:notation}
\chng{
We first recall some basic definitions and introduce our notation related to maps.
A map is a finite connected multigraph properly embedded in a surface that cuts the surface in a finite collection of simply connected domains, and is viewed up to orientation-preserving homeomorphisms of
the surface. 
We will think of each edge $\mve$ of a map as being a collection of two directed edges
$e$ and $\bar e$ that are paths in the underlying surface oriented in opposite directions, and we will write $\mve=\{e, \bar e\}$.
We will denote by $e(0)$ and $e(1)$ the origin (or initial vertex) and the target (or terminal vertex) of the directed edge $e$ respectively. 
Note that $e(0)=\bar e(1)$ and $e(1)=\bar e(0)$.
%We will picture the directed edge $e$ outgoing from $e(0)$ as being on the left of the  incoming directed edge $\bar e$.
For a map $\vm$, the set of all directed edges (resp. edges) of $\vm$ will be denoted by  
$\overrightarrow E (\vm)$ (resp. $E(\vm)$).
Thus, for $e\in\overrightarrow E (\vm)$, $\bar e$ will be the
corresponding element of $\overrightarrow E (\vm)$ that is oriented in
the opposite direction, and of course $\bar{\bar{e}}=e$. 
The image set of $e$ minus the points $e(0),e(1)$ is called the interior of
the edge $e$ (or of the corresponding unoriented edge $\mve$) and is denoted by $\INT(e)=\INT(\mve)$. 
A rooted map is a map $\vm$ together with a distinguished element $e_{\ast}\in \overrightarrow E (\vm)$. The latter is called the root edge, and $e_{\ast}(0)$ will be called the root vertex of $\vm$.
Combinatorially, the data of a map is equivalent to that of a finite connected multigraph where the set of directed edges emanating from each vertex is endowed with a cyclic order (that we think of as being clockwise). 	
For $e\in \overrightarrow E (\vm)$, we will write $e^-$ for the corner corresponding to $e$, that is a small angular sector between $f$ and $e$, where $f$ is the directed edge emanating from $e(0)$ that comes just before $e$ in the cyclic order around the vertex $e(0)$.
}

\chng{
We recall here that there is also an equivalent algebraic description in which a map $\vm$ is represented by a triple $(\alpha, \beta, \gamma)$ of permutations, where the cycles of $\alpha$, $\beta$, and $\gamma$ represent the edges, vertices, and faces of $\vm$ respectively. 
We refer the reader to \cite[Section~2]{chapuy} for a quick overview of this description.
We call $\gamma$ the face permutation of $\vm$.
In this representation of maps, a map is unicellular if and only if its face permutations consists of only one cycle.}

We make a note here that we use the notation $e, f$ to denote directed
edges of maps, and $\mve, \mvf$ to denote the corresponding edges. We
will also use the notation $e$ for a generic excursion, $f$ for a
generic function, and $\ve$ for a standard Brownian excursion. The
meaning will always be clear from the context, and there should not be
any confusion.

For $n\geq 1$ and $s\geq 0$, $\bM_{n, s}$ will denote the set of all rooted maps that have $n+1$ vertices and $n+s$ edges, where the root vertex has degree one. 
In particular, the root edge  cannot be a loop.
(The condition that the root vertex has degree one is of course artificial.
However, it will make certain things simpler.)
Thus, $\bM_{n, 0}$ is the set of plane trees on $n+1$ vertices, where the root vertex has degree one.

For a plane tree or a rooted labeled tree $\vt$, for every vertex $v$ of $\vt$, $\hght(v; \vt)$ will denote the tree distance between the root and $v$, and the vertices on the path connecting the root and $v$ (inclusive of both endpoints) will be called ancestors of $v$.
For $k=0, 1, \ldots$, $z(\vt; k)$ will denote the number of vertices in $\vt$ at height $k$.

\chng{We can explore any $\vt\in\bM_{n, 0}$ in the following fashion: 
Let $e_1$ be the root edge of $\vt$, and for $1\leq i\leq 2n-1$, let $e_{i+1}$ be the directed edge in $\vt$ with $e_{i+1}(0)=e_i(1)$ that comes right after $\bar e_i$ in the cyclic order on directed edges emanating from $e_i(1)$.
We say that the directed edges $e_1,\ldots, e_{2n}$ (resp. the corners $e_1^-,\ldots, e_{2n}^-$) appear in the contour order.
We refer to $e_1^-$ as the $0$-th corner of $\vt$.
Let $C_{\vt}(0)=0$ and $C_{\vt}(i)=\hght(e_i(1); \vt)$ for $1\leq i\leq 2n$.
Extend $C_{\vt}$ to a function on $[0, 2n]$ by linear interpolation.
Then $C_{\vt}:[0, 2n]\to[0,\infty)$ is called the contour function of $\vt$.}

For any set $A$, $\# A$ will denote its cardinality.
For two sequences $\{a_n\}_{n\geq 1}$ and $\{b_n\}_{n\geq 1}$ of positive numbers, we will write $a_n\sim b_n$ to mean that $a_n/b_n\to 1$ as $n\to\infty$.
Throughout this paper $c, c'$ will denote positive universal constants, and their values may change from line to line.

\section{Proofs}\label{sec:proofs}
The proofs of our main results will be given in this section.

\subsection{Exploration of maps}\label{sec:bf-df-search}
Fix integers $n, s\geq 1$, and $\vm\in \bM_{n, s}$. 
Let $e_{\ast}$ be the root edge of $\vm$, and let $\mve_{\ast}=\{e_{\ast}, \bar e_{\ast}\}$.

\medskip

\noindent{\bf Depth-first (DF) exploration:}
We will explore $\vm$ and simultaneously grow a plane tree $\vt$.
Set $\vt$ to be the edge $\mve_{\ast}$.
Root $\vt$ at $e_{\ast}$.
Set $e_1=e_{\ast}$.
Also, set the `current' map $\vm_{\cur}=\vm$.
Set $i=1$, and iterate as follows:

\begin{inparaenuma}
\noindent\item
If $i=2n$, stop.
Otherwise, \chng{let} $f$ be the directed edge \chng{in $\vm_{\cur}$ with $f(0)=e_i(1)$ that comes right after $\bar e_i$ in the cyclic order on directed edges emanating from $e_i(1)$},
and let $\mvf=\{f, \bar f\}$.  
\chng{(Note that if $\bar e_i$ is the only directed edge emanating from $e_i(1)$ in $\vm_{\cur}$, then $f=\bar e_i$.)}
Go to the next step. 

\noindent\item
If $\mvf$ is already an edge in $\vt$, set $e_{i+1}=f$, update $i$ to
$i+1$, and go to step (a). 
Otherwise, go to the next step. 

\noindent\item
If adding $\mvf$ to $\vt$ \chng{results in a tree}, then do so while respecting the circular order on directed edges emanating from $e_i(1)$ in $\vm$, set $e_{i+1}=f$, update $i$ to $i+1$, and go to step (a).
Otherwise, go to the next step.

\noindent\item
%If adding $\mvf$ to $\vt$ \chng{does not result in a tree}, 
Update $\vm_{\cur}$ to the map $\vm_{\cur}\setminus\INT(\mvf)$, and go to step (a).

\end{inparaenuma}

\medskip

After the algorithm terminates, we will get a plane tree $\vt\in\bM_{n, 0}$.
Note that $\vm_{\cur}=\vt$ at this stage.
In the contour exploration of $\vt$, its directed
edges will appear in the order $e_1,\ldots, e_{2n}$.
We can recover $\vm$ from $\vt$ by adding edges between certain corners in $\vt$.

Motivated by the last fact, we define the `depth-first admissible corners' of a plane tree as follows:
Fix a plane tree $\vt'\in\bM_{n, 0}$ and $s\geq 1$. 
Let $f_1,\ldots, f_{2n}$ be the directed edges of $\vt'$ in contour order.
Let $\dfac(\vt', s)$ be the set of all sequences $(i_1,\ldots, i_{2s}, k_1,\ldots,k_{2s})$, 
where $1\leq i_j\leq 2n-1$ and $1\leq k_j\leq 2s$ for $1\leq j\leq 2s$ such that the following hold:

\begin{itemize}
\item[{\bf (A.1)}]%\label{cond:1}
$i_{2j-1}\leq i_{2j}$ and $f_{i_{2j}}(0)$ is an ancestor of $f_{i_{2j-1}}(0)$, $1\leq j\leq s$.

\item[{\bf (A.2)}]%\label{cond:2}
For every maximal subset $\{r_1, \ldots, r_j\}\subseteq \{1,\ldots, 2s\}$ such that $i_{r_1}=\ldots=i_{r_j}$, $k_{r_1}, \ldots, k_{r_j}$ is a permutation of $1,\ldots, j$.

\item[{\bf (A.3)}]%\label{cond:3}
For $1\leq j<\ell\leq s$, either (a) $i_{2j-1}<i_{2\ell-1}$, or (b) $i_{2j-1}=i_{2\ell-1}$ and $i_{2j}<i_{2\ell}$, or (c) $i_{2j-1}=i_{2\ell-1}$ and $i_{2j}=i_{2\ell}$ and $k_{2j-1}<k_{2\ell-1}$.

\end{itemize}

The condition $i_{2j-1}\leq i_{2j}$ in {\bf (A.1)} means that the corner $f_{i_{2j}}^-$ appears after $f_{i_{2j-1}}^-$ in the contour order.
We would like to add edges between the corners \chng{$f_{i_{2j-1}}^-$ and $f_{i_{2j}}^-$}, $1\leq j\leq s$.
Condition {\bf (A.2)} gives a way of ordering the corresponding directed edges using the integers $k_j$ when multiple $i_j$-s are the same.
Condition {\bf (A.3)} is needed to pick out one particular representative among different possible permutations of the same sequence.

Given $\vt'\in\bM_{n, 0}$ and $\mvxi'=\big(i_j, k_j\, ;\ 1\leq j\leq
2s\big)\in\dfac(\vt', s)$, let $\cI(\vt', \mvxi')$ to be the map
obtained by adding edges between the corners \chng{$f_{i_{2j-1}}^-$ and $f_{i_{2j}}^-$}, $1\leq j\leq s$, while using the integers $k_j$, $1\leq j\leq 2s$, to order the corresponding directed edges when multiple $i_j$-s are the same.
For example, suppose $i_1=i_3$ and $k_3<k_1$. 
Let $\mvf_{12}=\{f_{12}, \bar f_{12}\}$ be the edge added between the corners $f_{i_1}^-$ and $f_{i_2}^-$ where $f_{12}$ is directed from $f_{i_1}^-$ towards $f_{i_2}^-$.
Similarly define $\mvf_{34}$.
Then in the resulting map, $f_{12}$ falls between $f_{34}$ and $f_{i_1}$ in the circular order on directed edges emanating from $f_{i_1}(0)$.
Similarly, if $i_1< i_3=i_2$ and $k_3<k_2$, then $\bar f_{12}$ falls between $f_{34}$ and $f_{i_2}$.
\chng{More generally, suppose $\{r_1, \ldots, r_j\}$ is a maximal subset of $\{1,\ldots, 2s\}$ such that $i_{r_1}=\ldots=i_{r_j}=i_{\star}$, say and  $\tilde r_1, \ldots, \tilde r_j$ is a permutation of $r_1,\ldots, r_j$ such that $k_{\tilde r_{\ell}}=\ell$ for $1\leq \ell\leq j$.
Then, in the cyclic order at $f_{i_{\star}}(0)$ in the map $\cI(\vt', \mvxi')$, the directed edge corresponding to $\tilde r_{\ell}$ is followed by the directed edge corresponding to $\tilde r_{\ell+1}$, $1\leq \ell\leq j-1$, and the directed edge corresponding to $\tilde r_j$ is followed by $f_{i_{\star}}$.}

\begin{figure}
	\centering
	\includegraphics[trim=1.7cm 15.5cm 3cm 2.2cm, clip=true, angle=0, scale=.5]{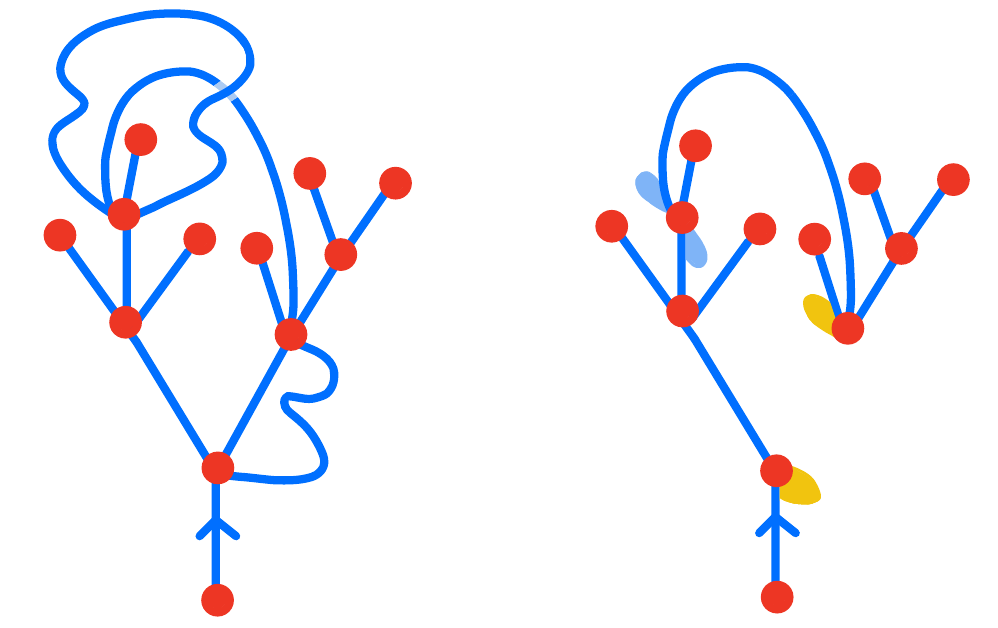}
	\captionof{figure}{On the left, a map $\vm$ with clockwise order at each vertex. On the right, $\df(\vm)$ with the relevant corners colored.}
	\label{fig:3-4}
\end{figure}

Now, if $\vt$ is the plane tree obtained from the depth-first exploration of $\vm\in\bM_{n, s}$ as above, then there exists a unique $\mvxi\in\dfac(\vt, s)$ such that $\cI(\vt, \mvxi)=\vm$. 
We set $\df(\vm)=(\vt,\mvxi)$.
An illustration is given in Figure \ref{fig:3-4}.
In this example, an edge is to be added between the two blue corners, and two edges are to be added between the yellow corners.
Here, $i_3=i_5$, $i_4=i_6$, $k_1=k_2=k_3=k_6=1$, and $k_4=k_5=2$.

\begin{rem}\label{rem:df-seaerch-variant}
In the above depth-first exploration, we follow the contour of the map (similar to the contour exploration of a plane tree), and delete any newly found edge if it completes a cycle in the map.
If $\vt$ is the plane tree resulting from this exploration, then each deleted edge is of the form $(u, v)$, where $v$ is an ancestor of $u$ in $\vt$.

There is a variant of the depth-first search where all the neighbors of the vertex being currently explored are immediately `discovered.'
If $\vt_1$ is the plane tree resulting from this exploration algorithm, then the edges deleted that are not loops will be of the form $(u, v)$, where $v$ is a child of some (strict) ancestor of $u$ in $\vt_1$, and further, $v$ lies to the right of the ancestral line of $u$.
To prove Theorem \ref{thm:bf=df-Hs} using this latter version of the depth-first search, we need some control over the number of vertices at distance $2$ from the ancestral line of a typical vertex in a uniform plane tree.
This can be done, but it will make the proof a bit more complicated.
\end{rem}

\noindent{\bf Breadth-first (BF) exploration:}
As before, we will explore $\vm$ and simultaneously grow a plane tree $\vt$.
Set $\vt$ to be the edge $\mve_{\ast}$ and root $\vt$ at $e_{\ast}$.
Set  $e_1=e_{\ast}$.
%Consider the directed edges with origin $e_{\ast}(0)$ that come after $e_{\ast}$ sequentially, and add the corresponding edges to $\vt$ if doing so does not create a cycle in $\vt$.
%Endow these edges with the circular order inherited from the order in $\vm$.
%Set $r$ to be the number of edges in $\vt$ at this stage.
%For $2\leq j\leq r$, set $e_j$ to be the $j$-th directed edge emenating from $e_{\ast}(0)$ in $\vt$.
Set $i=r=1$ and iterate as follows:

\begin{inparaenuma}
\noindent\item
If $r=n$, stop. 
Otherwise, go to the next step.
	
\noindent\item 
Consider the directed edges in $\vm$ with origin $e_i(1)$ that come after $\bar e_i$ sequentially.
Let $\Delta$ be the number of edges among them whose addition to $\vt$ results in a tree.
Update $\vt$ by adding these $\Delta$ edges to $\vt$ while maintaining the cyclic order inherited from the order in $\vm$.
For $1\leq j\leq\Delta$, set $e_{r+j}$ to be the $j$-th directed edge (after $\bar e_i$) emanating from $e_i(1)$ in $\vt$.
Update $i$ to $i+1$ and $r$ to $r+\Delta$.
Go to step (a).
\end{inparaenuma}

\medskip
After the algorithm terminates, we will get a plane tree $\vt\in\bM_{n, 0}$, and we can recover $\vm$ from $\vt$ by adding edges between certain corners in $\vt$.
Motivated by this, we define the `breadth-first admissible corners' of a plane tree as follows:
Fix a plane tree $\vt'\in\bM_{n, 0}$ and $s\geq 1$. 
Let $f_1,\ldots, f_{2n}$ be the directed edges of $\vt'$ in contour order.
Let $\bfac(\vt', s)$ be the set of all sequences $(i_1,\ldots, i_{2s}, k_1,\ldots,k_{2s})$, where $1\leq i_j\leq 2n-1$ and $1\leq k_j\leq 2s$ for $1\leq j\leq 2s$ such that Conditions {\bf (A.2)} and {\bf (A.3)} hold as above, and further the following holds:
For $1\leq j\leq s$, $i_{2j-1}\leq i_{2j}$ and 
	\begin{align}\label{eqn:16}
\hght\big(f_{i_{2j}}(0)\big)\in
\big\{ \hght\big(f_{i_{2j-1}}(0)\big)\, ,\ \hght\big(f_{i_{2j-1}}(0)\big)-1 \big\}\, .
	\end{align}

The only difference from the depth-first case is \eqref{eqn:16}.
As before, we would like to add edges between the corners \chng{$f_{i_{2j-1}}^-$ and $f_{i_{2j}}^-$}, $1\leq j\leq s$.
Thus, \eqref{eqn:16} essentially says that edges are to be added between vertices at roughly the same height.

Given $\vt'\in\bM_{n, 0}$ and $\mvxi'=\big(i_j, k_j\, ;\ 1\leq j\leq 2s\big)\in\bfac(\vt', s)$, let $\cI(\vt', \mvxi')$ to be the map obtained by adding edges between the corners \chng{$f_{i_{2j-1}}^-$ and $f_{i_{2j}}^-$}, $1\leq j\leq s$, and as before, we use the integers $k_j$, $1\leq j\leq 2s$, to order the corresponding directed edges when multiple $i_j$-s are the same.

\begin{figure}
	\centering
	\includegraphics[trim=1.7cm 16cm 3cm 1.75cm, clip=true, angle=0, scale=.5]{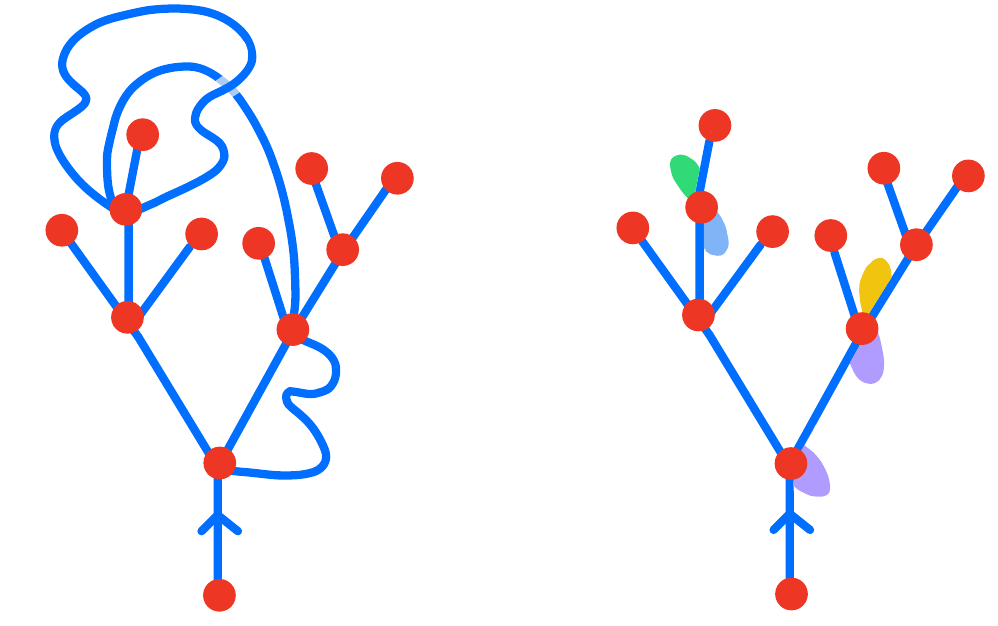}
	\captionof{figure}{On the left, a map $\vm$ with clockwise order at each vertex. On the right, $\bbf(\vm)$ with the relevant corners colored.}
	\label{fig:5-6}
\end{figure}

Now, if $\vt$ is the plane tree obtained from the breadth-first exploration of $\vm\in\bM_{n, s}$, then there exists a unique $\mvxi\in\bfac(\vt, s)$ such that $\cI(\vt, \mvxi)=\vm$. 
We set $\bbf(\vm)=(\vt,\mvxi)$.
An illustration is given in Figure \ref{fig:5-6}. 
In this example, an edge is to be added between the green and the blue corner, the green and the yellow corner, and the two purple corners.
Here, $i_1=i_3$ correspond to the green corner,  $i_2$ and $i_4$ correspond to the blue and yellow corners respectively, $i_5$ and $i_6$ correspond to the purple corners, and
$k_1=1$, $k_3=2$, $k_2=k_4=k_5=k_6=1$.

We define 
\begin{gather}\label{eqn:17}
\bft(n,s):=\big\{ \big(\vt, \mvxi\big)\, :\, \vt\in\bM_{n, 0}\, ,\ \mvxi\in\bfac(\vt, s) \big\}\, ,\ \text{ and}\\
\dft(n,s):=\big\{ \big(\vt, \mvxi\big)\, :\, \vt\in\bM_{n, 0}\, ,\ \mvxi\in\dfac(\vt, s) \big\}\, .
\end{gather}

\subsection{Proof of Theorem \ref{thm:bf=df-Hs}}\label{sec:proof-bf=df-Hs}
The proof relies on the following two propositions.

\begin{prop}\label{prop:bf-scaling-limit}
	Fix an integer $s\geq 1$.
	For each $n\geq 1$, let $\big(T_{n}^{\bbf}\, ,\,  \Xi_{n}^{\bbf}\big)$ be uniformly distributed over $\bft(n,s)$.
	Let $G_{n, s}^{\bbf}$ be the rooted metric measure space obtained by endowing 
	$\cI\big(T_{n}^{\bbf}\, ,\,  \Xi_{n}^{\bbf}\big)$ with the graph distance, declaring the root vertex as the root, and finally assigning probability $1/n$ to every vertex of $\cI\big(T_{n}^{\bbf}\, ,\,  \Xi_{n}^{\bbf}\big)$ except the root vertex.
	Then
	\[\frac{1}{\sqrt{n}}\cdot  G_{n, s}^{\bbf}\weakc  \frac{1}{\sqrt{2}}\cdot\hsbf\, ,\ \ \text{ as }\ \ n\to\infty\]
	w.r.t. the pointed GHP topology.
\end{prop}

\begin{prop}\label{prop:df-scaling-limit}
	Fix an integer $s\geq 1$.
	For each $n\geq 1$, let $\big(T_{n}^{\df}\, ,\,  \Xi_{n}^{\df}\big)$ be uniformly distributed over $\dft(n,s)$.
	Let $G_{n, s}^{\df}$ be the rooted metric measure space obtained by endowing 
	$\cI\big(T_{n}^{\df}\, ,\,  \Xi_{n}^{\df}\big)$ with the graph distance, declaring the root vertex as the root, and finally assigning probability $1/n$ to every vertex of $\cI\big(T_{n}^{\df}\, ,\,  \Xi_{n}^{\df}\big)$ except the root vertex.
	Then
	\[\frac{1}{\sqrt{n}}\cdot G_{n, s}^{\df}\weakc \frac{1}{\sqrt{2}}\cdot\cH_{(s)}\, ,\ \ \text{ as }\ \ n\to\infty\]
	w.r.t. the pointed GHP topology.
\end{prop}

We first prove Theorem \ref{thm:bf=df-Hs} using Propositions \ref{prop:bf-scaling-limit} and \ref{prop:df-scaling-limit}.
The proofs of these propositions will be given in the Sections \ref{sec:proof-bf-scaling-limit} and \ref{sec:proof-df-scaling-limit} respectively.

\medskip

\noindent{\bf Completing the proof of Theorem \ref{thm:bf=df-Hs}:}
Fix $s\geq 1$ as the result is trivial for $s=0$.
Let $M_{n, s}\sim\mathrm{Uniform}(\bM_{n, s})$.
Since $\df: \bM_{n, s}\to\dft(n, s)$ and $\bbf: \bM_{n, s}\to\bft(n, s)$ are bijections with inverse $\cI$, 
\begin{align}\label{eqn:18}
\cI\big(T_{n}^{\df}\, ,\,  \Xi_{n}^{\df}\big)
\equald M_{n, s}
\equald
\cI\big(T_{n}^{\bbf}\, ,\,  \Xi_{n}^{\bbf}\big)
\, .
\end{align}
View $M_{n, s}$ as a rooted metric measure space by endowing it with
the graph distance, declaring the root vertex as the root, and assigning probability $1/n$ to every non-root vertex.
Then \eqref{eqn:18} implies that $G_{n, s}^{\df}\equald M_{n, s}\equald G_{n, s}^{\bbf}$.
Thus using Propositions \ref{prop:bf-scaling-limit} and \ref{prop:df-scaling-limit}, we get
\begin{align}\label{eqn:92}
\sqrt{\frac{2}{n}}\cdot M_{n, s}\weakc\cH_{(s)}\equald\hsbf\, 
\end{align}
w.r.t. pointed GHP topology, which completes the proof.
\qed

\medskip

\begin{rem}\label{rem:uniform-map-scaling-limit}
Let $M_{n, s}'$ be uniformly distributed over the set of all rooted maps on $n$ vertices having $s$ surplus edges (i.e., we drop the condition that the degree of the root vertex is one).
Then clearly, \eqref{eqn:92} continues to hold if we replace $M_{n, s}$ by $M_{n, s}'$.
\end{rem}

\subsection{Proof of Proposition \ref{prop:bf-scaling-limit}}\label{sec:proof-bf-scaling-limit}
\begin{figure}
	\centering
		\includegraphics[trim=1.7cm 19cm 3cm 2.4cm, clip=true, angle=0, scale=.8]{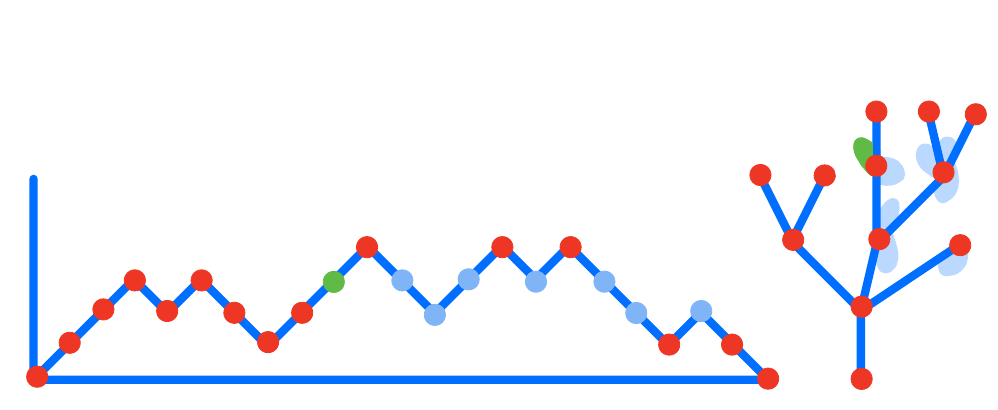}
		\captionof{figure}{Left: $f$ with the point on the graph of $f$ corresponding to $i=9$ colored green, and points that correspond to $\fB(f; 9)\setminus\{9\}$ colored light blue. 
		Right: the corresponding plane tree with the relevant corners colored.}
		\label{fig:7}
\end{figure}

Throughout this section we work with a fixed $s\geq 1$.
For $n\geq 1$, define
\begin{align}
\fC_{n}:=\big\{f:[0, 2n]\to [0,\infty) \ 
:\ & f(0)=f(2n)=0\, ,\ f(i)\in\bZ_{>0} \text{ for }  1\leq i\leq 2n-1\, ,\notag \\
& |f(i+1)-f(i)|=1 \text{ for }  0\leq i\leq 2n-1\, ,\text{ and } f(t) \text{ is given}\notag\\
&\text{by linear interpolation for other values of } t\in[0, 2n] \big\}.\label{eqn:52}
\end{align}
Note that $\fC_{n}$ is the set of all contour functions of plane trees in $\bM_{n, 0}$.
For $f\in\fC_{n}$ and $i=0, 1,\ldots, 2n$, let
\begin{gather*}
\fB(f; i):= \big\{j\in\{i\vee 1,\ldots, 2n-1\}\, :\, f(j)=f(i)\, \text{ or }\, f(i)-1\big\}\, ,\\
B(f; i):=\#\, \fB(f; i)\, , \ \text{ and }\ B(f)=\sum_{i=0}^{2n}B(f, i)\, .
\end{gather*}
An illustration of $\fB(f; i)$ is given in Figure \ref{fig:7}.
%In terms of the plane tree whose contour process is $f$, for $1\leq i\leq 2n-1$, $\fB(f; i)\setminus\{i\}$ corresponds to the set of corners that are at height $f(i)$ or $f(i)-1$ and appear strictly after the $i$-th corner in the contour order.
Note the connection with the set of breadth-first admissible corners $\bfac(\cdot, \cdot)$.

Let $C_n^{\circ}$ be distributed as
\begin{align}\label{eqn:19}
\pr\big(C_{n}^{\circ}=f\big)
=
\frac{\big(B(f)\big)^s}{\sum_{\phi\in\fC_{n}}\big(B(\phi)\big)^s}\, , \ \ \ \  f\in\fC_{n}\, .
\end{align}
Conditional on $C_n^{\circ}$, sample i.i.d. random variables 
$\big(I_{n, 2j-1}^{\circ}\, ,\, I_{n,2j}^{\circ}\big)$, $1\leq j\leq s$, where
\begin{align}\label{eqn:20}
\pr\big(I_{n, 1}^{\circ}=i\big)
=
B(C_n^{\circ}; i)/B(C_n^{\circ})\, , \ \ i=0, 1,\ldots, 2n\, ,
\end{align}
and conditional on $I_{n, 1}^{\circ}$, $I_{n, 2}^{\circ}$ is uniform over $\fB(C_n^{\circ}; I_{n, 1}^{\circ})$.

Let $T_n^{\circ}$ be the plane tree whose contour function is $C_n^{\circ}$.
For $1\leq j\leq s$, add an edge to $T_n^{\circ}$ between the $I_{n, 2j-1}^{\circ}$-th corner and the $I_{n,2j}^{\circ}$-th corner of $T_n^{\circ}$.
Endow the resulting space \greg{with} the graph distance, root it at the root vertex of $T_n^{\circ}$, and assign probability $1/n$ to every vertex of $T_n^{\circ}$ except the root vertex. 
Denote the resulting rooted metric measure space by $G_{n, s}^{\circ}$.
We complete the proof of Proposition \ref{prop:bf-scaling-limit} by combining the next two lemmas:

\begin{lem}\label{lem:1}
We have, 
\[
\frac{1}{\sqrt{n}}\cdot G_{n, s}^{\circ}\weakc \frac{1}{\sqrt{2}}\cdot\cH_{(s)}^{\bbf},\ \ \text{ as }\ \ n\to\infty
\]
w.r.t. the pointed GHP topology.
\end{lem}
\begin{lem}\label{lem:2}
There exists a coupling of $G_{n, s}^{\circ}$ and $G_{n, s}^{\bbf}$ such that
\[
\pr\big(G_{n, s}^{\circ}\neq G_{n, s}^{\bbf}\big)\to  0,\ \ \text{ as }\ \ n\to\infty\, .
\]
\end{lem}
The rest of this section is devoted to the proofs of Lemmas \ref{lem:1} and \ref{lem:2}.
To this end, let us define the `discrete local time' of $f\in\fC_n$ as
\begin{align}\label{eqn:21}
L(f; t, y):=\#\big\{0\leq j\leq t\, :\, f(j)=y\big\}\, ,\ \ t=0, 1,\ldots, 2n\, , \ \ y\in\bZ\, ,
\end{align}
and extend it to a function on $[0, 2n]\times\bR$ via the relation
\begin{align}\label{eqn:22}
L(f; t, y):=
&  
\big(\langle t\rangle -t \big) \big(\langle y\rangle-y \big) L \big(f; \lfloor t\rfloor, \lfloor y\rfloor \big) 
+ 
\big(\langle  t \rangle -t\big) \big(y-\lfloor y\rfloor\big) L\big(f; \lfloor t\rfloor, \langle y\rangle\big)\notag 
\\
&
+
\big(t-\lfloor t\rfloor\big) \big(\langle y\rangle-y \big) L \big(f; \langle t\rangle, \lfloor y\rfloor \big) 
+ 
\big(t-\lfloor t\rfloor\big) \big(y-\lfloor y\rfloor\big) L\big(f; \langle t\rangle, \langle y\rangle\big)
\, ,
\end{align}
where $\lfloor t\rfloor:=\max\big\{j\in\bZ\, :\, j\leq t\big\}$ and $\langle t\rangle :=\lfloor t\rfloor+1$.
Note that $L(f; \cdot,\cdot)$ is a continuous function on $[0, 2n]\times\bR$.

Let $C_n$ be uniformly distributed over $\fC_n$.
Define
\begin{gather}
L_n(\cdot, \cdot):=L(C_n;\, \cdot ,\, \cdot)\, ,\ \ 
\| L_n\|_{\infty}:=\sup_{y\in\bR} L_n(2n, y), \ \ \text{ and }\label{eqn:26}\\
\bar L_n(t, y):=(2n)^{-1/2} L_n\big(2nt,\, y\sqrt{2n}\big)\, ,\ 0\leq t\leq 1\, ,\, y\in\bR\, .\label{eqn:27}
\end{gather}
Further, let
\begin{align}\label{eqn:23}
\bar C_n(t):=(2n)^{-1/2} C_n(2nt)\, ,\ 0\leq t\leq 1\, .
\end{align}
We will make use of the following result in our proof.
\begin{prop}\label{prop:local-time-convergence}
The following convergence holds in $C[0, 1]\times C\big([0, 1]\times\bR\big)$:
\begin{gather}
\big(\bar C_n,\, \bar L_n \big)
\weakc 
\big(\ve(\cdot),\, 
\eta(\ve;  \cdot, \cdot)\big)\, .\label{eqn:24}
\end{gather}
Further, there exist constants $c_1, c_2 >0$ such that for all $n\geq 1$,
\begin{align}\label{eqn:25}
\pr\big(\|L_n\|_{\infty}\geq u\sqrt{n}\big)\leq c_1\exp(-c_2 u^2)\, , \ \  u\geq 0\, .
\end{align}
\end{prop}
Let $T_n$ be the plane tree whose contour function is $C_n$, and let $e_{\ast}$ be its root edge.
Then the tree obtained from $T_n$ by deleting $e_{\ast}(0)$ and $\INT(e_{\ast})$ has the same distribution as a uniform plane tree on $n$ vertices.
The convergence $\bar C_n\weakc\ve$ can be deduced from the fact that simple random walk excursions converge, after proper rescaling, to $\ve$, 
or by using the relation between $T_n$ and uniform plane trees and \cite[Theorem 23]{aldous-crt-3}.
We now state a general result that implies \eqref{eqn:25}. 
This result will also be used in later sections.
Recall the notation $z(\vt; k)$ from Section \ref{sec:notation}.

\begin{thm}[\cite{addarioberry-devroye-janson}, Theorem 1.1]\label{thm:conditioned-gw-tree}
Let $\xi$ be a nonnegative integer-valued random variable with $\bE\xi=1$ and $0<\var\xi<\infty$.
Let $T_{\xi}$ be a plane Galton-Watson tree with offspring distribution $\xi$.
Let $|T_{\xi}|$ denote the number of vertices in $T_{\xi}$.
Then there exist $c_1, c_2>0$ such that 
\[
\pr\big(\max_k z(T_{\xi}; k)\geq x\sqrt{n}\ \big|\ |T_{\xi}|=n\big)\leq c_1\exp\big(-c_2 x^2\big)\, ,
\]
for all $x>0$ and $n\geq 1$ such that $\pr\big(|T_{\xi}|=n\big)>0$.
\end{thm}
Using the connection between $T_n$ and uniform plane trees, \eqref{eqn:25} follows if we apply Theorem \ref{thm:conditioned-gw-tree} with $\mathrm{Geometric}(1/2)$ offspring distribution and the fact that $L_n(2n, k)= z(T_n, k)+z(T_n, k+1)$ for $k\geq 0$.
A brief sketch of the proof of the full convergence in \eqref{eqn:24} is given in Appendix \ref{sec:appendix}.
Let us now prove Lemma \ref{lem:1}.

\medskip

\noindent {\bf Proof of Lemma \ref{lem:1}:}
Recall the definition of $C_n^{\circ}$ from \eqref{eqn:19}, and 
analogous to \eqref{eqn:23}, \eqref{eqn:26} and \eqref{eqn:27}, define $\bar C_n^{\circ}$, $L_n^{\circ}$, and $\bar L_n^{\circ}$.
Note that
\begin{align}\label{eqn:28}
\bE\big[\phi\big(\bar C_n^{\circ}\, ,\, \bar L_n^{\circ}\big)\big]
=
\frac{\bE\big[\phi\big(\bar C_n\, ,\, \bar L_n\big)B(C_n)^s\big]}{\bE\big[B(C_n)^s\big]}
\end{align}
for any bounded continuous $\phi:C[0, 1]\times C\big([0, 1]\times\bR\big)\to\bR$.
Now
\begin{align}\label{eqn:30}
B(C_n)=
\sum_{i=1}^{2n-1}
\big[
L_n\big(2n, C_n(i)\big)-L_n\big(i-1, C_n(i)\big)
+
L_n\big(2n, C_n(i)-1\big)-L_n\big(i-1, C_n(i)-1\big)
\big]\, .
\end{align}
Thus, the following convergence holds jointly with \eqref{eqn:24}:
\begin{align}\label{eqn:29}
\frac{B(C_n)}{(2n)^{3/2}}\weakc
& \ 2\int_{0}^{1}\big[\eta(\ve; 1, \ve(t))-\eta(\ve; t, \ve(t))\big]dt\notag\\
& =2\int_{y\in\bR}dy\int_{t=0}^{1}\big[\eta(\ve; 1, y)-\eta(\ve; t, y)\big]\eta(\ve; dt, y)
=\int_{y\in\bR}\eta(\ve; 1, y)^2dy\, ,
\end{align}
where the convergence in the first step uses \eqref{eqn:24}, the second step uses \eqref{eqn:6}, and the last step follows from a direct computation.
Now, \eqref{eqn:30} implies that
$
B(C_n)\leq 4n\|L_n\|_{\infty}
$.
Combined with \eqref{eqn:25}, this shows that for any $s\geq 1$, the sequence of random variables 
\begin{align}\label{eqn:44}
\big\{(2n)^{-3s/2}B(C_n)^s\big\}_{n\geq 1}\ \text{ is uniformly integrable.}
\end{align}
Hence, we conclude from \eqref{eqn:28}, \eqref{eqn:29}, and \eqref{eqn:24} that
\begin{align*}
\bE\big[\phi\big(\bar C_n^{\circ}\, ,\, \bar L_n^{\circ}\big)\big]
\weakc
\bE\big[\phi\big(\esbf(\cdot)\, ,\, \nsbf(\cdot, \cdot)\big)\big]\, ,
\end{align*}
where $\nsbf$ is as in \eqref{eqn:38}.
From now on, we work in a space where the following almost sure convergence holds in $C[0, 1]\times C\big([0, 1]\times\bR\big)$:
\begin{align}\label{eqn:31}
\big(\bar C_n^{\circ}\, ,\, \bar L_n^{\circ}\big)
\convas
\big(\esbf\, ,\, \nsbf\big)\, .
\end{align}
Similar to \eqref{eqn:29},
using \eqref{eqn:31}, we see that 
\begin{align}\label{eqn:32}
(2n)^{-3/2}B(C_n^{\circ})
\convas
\int_0^{\infty}\nsbf(1, y)^2dy\, .
\end{align}
From \eqref{eqn:20}, the law of $I_{n, 1}^{\circ}$ is given by
\[
\pr\big(I_{n, 1}^{\circ}\geq i\ \big|\ C_n^{\circ} \big)=\frac{1}{B(C_n^{\circ})}\sum_{j=i}^{2n}
\big[
L_n^{\circ}\big(2n, C_n^{\circ}(j)\big)-L_n^{\circ}\big(j-1, C_n^{\circ}(j)\big)
+
L_n^{\circ}\big(2n, C_n^{\circ}(j)-1\big)-L_n^{\circ}\big(j-1, C_n^{\circ}(j)-1\big)
\big]
\]
for $i=0, 1,\ldots, 2n$.
Using \eqref{eqn:31} and \eqref{eqn:32}, we see that $(2n)^{-1}I_{n, 1}^{\circ}$ converges in distribution to a random variable $u_{(1)}$ with law given by
\begin{align}\label{eqn:35}
&\pr\big(u_{(1)}\geq u\ \big|\ \esbf\big)\cdot \int_0^{\infty}\nsbf(1, y)^2dy 
=
2\cdot
\int_{t=u}^1\Big[\nsbf\big(1,\, \esbf(t)\big)-\nsbf\big( t,\, \esbf(t)\big)\Big]dt
\notag\\
&\hskip50pt
=
2\cdot
\int_{y=0}^{\infty}dy\int_{t=u}^{1}
\Big[\nsbf\big(1, y\big)-\nsbf\big(t, y\big)\Big]\nsbf\big( dt, y\big)
\notag\\
&\hskip100pt
=\int_{y=0}^{\infty}\Big[\nsbf\big(1, y\big)-\nsbf\big(u, y\big)\Big]^2 dy \, ,
\ \ \ 0\leq u\leq 1\, ,
\end{align}
where the second equality uses \eqref{eqn:6}, and the last step follows from a direct computation.
We can assume that we are working in a space where in addition to \eqref{eqn:31},
\begin{align}\label{eqn:33}
(2n)^{-1}I_{n, 1}^{\circ}\convas u_{(1)}\,  .
\end{align}
Conditional on $I_{n, 1}^{\circ}$, 
$I_{n, 2}^{\circ}$ is uniformly distributed over $\fB(C_n^{\circ}; I_{n, 1}^{\circ})$.
Thus, the conditional distribution function of $I_{n, 2}^{\circ}$ is given by
\[
\frac{
L_n^{\circ}\big(j, C_n^{\circ}(I_{n, 1}^{\circ})\big)
-L_n^{\circ}\big(I_{n, 1}^{\circ}-1\, ,\, C_n^{\circ}(I_{n, 1}^{\circ})\big)
+L_n^{\circ}\big(j, C_n^{\circ}(I_{n, 1}^{\circ})-1\big)
-L_n^{\circ}\big(I_{n, 1}^{\circ}\, ,\, C_n^{\circ}(I_{n, 1}^{\circ})-1\big)
}{
L_n^{\circ}\big(2n, C_n^{\circ}(I_{n, 1}^{\circ})\big)
-L_n^{\circ}\big(I_{n, 1}^{\circ}-1\, ,\, C_n^{\circ}(I_{n, 1}^{\circ})\big)
+L_n^{\circ}\big(2n, C_n^{\circ}(I_{n, 1}^{\circ})-1\big)
-L_n^{\circ}\big(I_{n, 1}^{\circ}\, ,\, C_n^{\circ}(I_{n, 1}^{\circ})-1\big)
}
\, ,
\]
for $j=I_{n,1}^{\circ}, I_{n,1}^{\circ}+1,\ldots, 2n$, and is a right continuous step function in between.
Thus, using \eqref{eqn:31} and \eqref{eqn:33}, we see that $(2n)^{-1}I_{n, 2}^{\circ}$ converges in distribution to a random variable $v_{(1)}$ with conditional distribution function (given $\esbf$ and $u_{(1)}$) given by
\begin{align}\label{eqn:36}
\frac{
\nsbf\big(v,\, \esbf(u_{(1)})\big)-\nsbf\big(u_{(1)},\, \esbf(u_{(1)})\big)
}{
\nsbf\big(1,\, \esbf(u_{(1)})\big)-\nsbf\big(u_{(1)},\, \esbf(u_{(1)})\big)
}\, ,\ \ 
u_{(1)}\leq v\leq 1\, .
\end{align}
Repeating the argument, we can assume that we are working in a space where
\begin{align}\label{eqn:34}
(2n)^{-1}\big(I_{n,1}^{\circ}, I_{n, 2}^{\circ}, \ldots, I_{n,2s-1}^{\circ}, I_{n, 2s}^{\circ}\big)
\convas
\big(u_{(1)}, v_{(1)}, \ldots, u_{(s)}, v_{(s)}\big)
\end{align}
in addition to \eqref{eqn:31}.
Here,  $(u_{(i)}, v_{(i)})$, $1\leq i\leq s$, are i.i.d. conditionally on $\esbf$.

It follows from \eqref{eqn:35} and \eqref{eqn:36} that conditional on $\esbf$,
\begin{align}\label{eqn:37}
\big(u_{(i)},\, v_{(i)}\big)
\equald
\big(\min\{u_i, v_i\},\, \max\{u_i, v_i\}\big)\, ,\ \ 1\leq i\leq s\, ,
\end{align}
where $u_i, v_i$, $1\leq i\leq s$, are as in Construction \ref{constr:H-s-bf}.
In particular, consider the rooted metric measure space obtained by identifying $q_{\esbf}(u_{(i)})$ and $q_{\esbf}(v_{(i)})$, $1\leq i\leq s$, in $\cT_{\esbf}$.
This space has the same distribution as $1/2\cdot\hsbf$.
We denote this space as $1/2\cdot\hsbf$ for the rest of this proof.
Let $\psi:[0, 1]\to 1/2\cdot\hsbf$ denote the quotient map.

Let $\bar G_{n, s}$ be the rooted metric measure space obtained from $\cT_{\bar C_n^{\circ}}$ by identifying $q_{\bar C_n^{\circ}}\big(I_{n,2j-1}^{\circ}/(2n)\big)$ and 
 $q_{\bar C_n^{\circ}}\big(I_{n,2j}^{\circ}/(2n)\big)$, $1\leq j\leq s$.
Let $\psi_n:[0, 1]\to\bar G_{n, s}$ denote the quotient map.

Let $\pi_n$ be the measure on $(1/2\cdot\hsbf)\times\bar G_{n,s}$ given by the push-forward of the Lebesgue measure on $[0, 1]$ under $\psi\times\psi_n$, 
and let $\cR_n$ be the correspondence \chng{(see, e.g., \cite[Section~2.1]{AddBroGolMie13} or \cite[Section~7.3.3]{burago-burago-ivanov} for the definitions of a correspondence and its distortion)} between  $1/2\cdot\hsbf$ and 
$\bar G_{n,s}$ given by 
\[
\cR_n:=\big\{ (x, y)\ :\ \exists\ t\in[0, 1]\text{ with }\psi(t)=x\text{ and }\psi_n(t)=y  \big\}\, .
\]
Then using \eqref{eqn:34} and the convergence $\bar C_n^{\circ}\convas\esbf$ from \eqref{eqn:31}, it is easy to see that the distortion of $\cR_n$ satisfies
$
\dis(\cR_n)\convas 0\, .
$
Further, $\pi_n(\cR_n^c)=0$, and the projections of $\pi_n$ onto $1/2\cdot\hsbf$ and $\bar G_{n, s}$ are same as the measures on the respective spaces.
Combining these observations, it follows that
\[
d_{\GHP}^1\big(\bar G_{n, s}\, ,\, 1/2\cdot\hsbf\big)\convas 0\, .
\]
There is an obvious coupling between $G_{n, s}^{\circ}$ and $\bar G_{n, s}$, and the proof will be complete if we show that
\begin{align}\label{eqn:39}
d_{\GHP}^1\big((2n)^{-1/2}\cdot G_{n, s}^{\circ}\, ,\, \bar G_{n, s}\big)\convas 0
\end{align}
in this coupling.
There are two differences between $(2n)^{-1/2}\cdot G_{n, s}^{\circ}$ and $\bar G_{n, s}$: 
(i) The measure on $\bar G_{n, s}$ is the normalized line measure, whereas the measure on $(2n)^{-1/2}\cdot G_{n, s}^{\circ}$ is the uniform probability measure on all non-root vertices of $T_n^{\circ}$.
(ii) The points $q_{\bar C_n^{\circ}}\big(I_{n,2j-1}^{\circ}/(2n)\big)$ and 
$q_{\bar C_n^{\circ}}\big(I_{n,2j}^{\circ}/(2n)\big)$, $1\leq j\leq s$, are identified in the construction of $\bar G_{n,s}$, whereas an edge of length $(2n)^{-1/2}$ is added between them in the construction of $(2n)^{-1/2}G_{n,s}^{\circ}$.
Using these observations, the proof of \eqref{eqn:39} is routine.
We omit the details.
\qed

\vskip5pt

Recall that $T_n$ denotes the plane tree whose contour function is $C_n$.
The proof of Lemma \ref{lem:2} relies on the following result:

\begin{lem}\label{lem:3}
We have,
\[
n^{-3s/2}\cdot\bE\big|s!\times\#\bfac(T_n, s)- B(C_n)^s\big|\to 0\ \ \text{ as }\ \ n\to\infty\, .
\]
\end{lem}

\noindent{\bf Proof:} Note that
\begin{align}\label{eqn:40}
\bfac(T_n, s)=\big\{&\big(i_1,\ldots, i_{2s}, 1,\ldots, 1\big)\ :\ 1\leq i_1<i_3<\ldots<i_{2s-1}\leq 2n-1\,
,\ \ \ 
i_{2j}\in\fB\big(C_n; i_{2j-1}\big)\notag\\
&\hskip100pt 
\text{ for }1\leq j\leq s\, , \text{ and } i_{\ell}\text{-s are all distinct}
\big\}\notag\\
&\bigcup
\big\{\big(i_1,\ldots, i_{2s}, k_1,\ldots, k_{2s}\big)\in\bfac(T_n, s)\ :\  i_{\ell}\text{-s are not all distinct}
\big\}\, .
\end{align}
Similarly,
\begin{align*}
B(C_n)^s=&\#\big\{(i_1,\ldots, i_{2s})\in[2n]^{2s}\ :\  
i_{2j}\in\fB\big(C_n; i_{2j-1}\big)
\text{ for }1\leq j\leq s\, , \text{ and } i_{\ell}\text{-s are all distinct}
\big\}\\
&+
\#\big\{(i_1,\ldots, i_{2s})\in[2n]^{2s}\ :\  
i_{2j}\in\fB\big(C_n; i_{2j-1}\big)
\text{ for }1\leq j\leq s\, , \text{ and } i_{\ell}\text{-s are not all distinct}
\big\}\, .
\end{align*}
Abbreviating `are not all distinct' as `NAD' 
and writing $(i_1,\ldots, i_{2s})$ for $(i_1,\ldots, i_{2s})\in [2n]^{2s}$,
we see that
\begin{align}\label{eqn:41}
&
0
\leq
s!\times\#\bfac(T_n, s)-B(C_n)^s
\leq 
s!\cdot\#\big\{\big(i_1,\ldots, i_{2s}, k_1,\ldots, k_{2s}\big)\in\bfac(T_n, s)\ :\  i_{\ell}\text{-s NAD}
\big\}\notag\\
&\hskip60pt
\leq
s!\times(2s)^{2s}\times
\#\big\{(i_1,\ldots, i_{2s})\ :\  
i_{2j}\in\fB\big(C_n; i_{2j-1}\big)\text{ for } 1\leq j\leq s\, , \text{ and }
i_{\ell}\text{-s NAD}
\big\}\notag\\
&\hskip120pt
\leq
s!\times(2s)^{2s}\times\big(A_{12}+A_{13}+A_{14}+A_{24}\big)\, ,
\end{align}
where
\begin{gather*}
A_{12}=\#\big\{ 
(i_1,\ldots, i_{2s})\ :\  
i_{2j}\in\fB\big(C_n; i_{2j-1}\big)\text{ for } 1\leq j\leq s\, , \
i_{2\ell-1}=i_{2\ell}\text{ for some }\ell
\big\}\, ,\\
A_{13}=\#\big\{ 
(i_1,\ldots, i_{2s})\ :\  
i_{2j}\in\fB\big(C_n; i_{2j-1}\big)\text{ for } 1\leq j\leq s\, , \
i_{2\ell-1}=i_{2\ell'-1}\text{ for some }\ell\neq\ell'
\big\}\, ,
\end{gather*}
and $A_{14}$ and $A_{24}$ are defined similarly with the respective defining conditions being 
`$i_{2\ell-1}=i_{2\ell'}$ for some $\ell\neq\ell'$' and
`$i_{2\ell}=i_{2\ell'}$ for some $\ell\neq\ell'$.'
Now,
\begin{align*}
A_{12}
&
\leq
s\cdot\#\big\{ 
(i_1,\ldots, i_{2s})\ :\  
i_{2j}\in\fB\big(C_n; i_{2j-1}\big)\text{ for } 1\leq j\leq s\, , \
i_1=i_2
\big\}\\
&
\leq 
s\cdot (2n)\cdot\big[2n\cdot\max_i B(C_n; i) \big]^{s-1}
\leq
s\cdot (2n)\cdot\big[2n\cdot 2\|L_n\|_{\infty}\big]^{s-1}\, ,
\end{align*}
where the last step uses the fact
\begin{align}\label{eqn:43}
\max_i B(C_n; i)\leq 2\|L_n\|_{\infty}\, .
\end{align}
Combined with \eqref{eqn:25}, we get $n^{-3s/2}\cdot\bE[A_{12}]=O(n^{-1/2})$.
Similarly,
\begin{align}\label{eqn:46}
A_{13}
&
\leq
s^2\cdot\#\big\{ 
(i_1,\ldots, i_{2s})\ :\  
i_{2j}\in\fB\big(C_n; i_{2j-1}\big)\text{ for } 1\leq j\leq s\, , \
i_1=i_3
\big\}\notag\\
&
\leq
s^2\cdot\#\big\{(i, i_2, i, i_4)\ :\  i_2,  i_4\in\fB\big(C_n; i\big)\big\}\times
\big[2n\cdot 2\|L_n\|_{\infty} \big]^{s-2}
\notag\\
&
\leq 
s^2\cdot (2n)\cdot\big(2\cdot\|L_n\|_{\infty}\big)^2\times
\big[2n\cdot 2\|L_n\|_{\infty} \big]^{s-2}\, .
\end{align}
Using \eqref{eqn:25} again, we get $n^{-3s/2}\cdot\bE[A_{13}]=O(n^{-1})$.
We can similarly show that $n^{-3s/2}\cdot\bE[A_{14}+A_{24}]=O(n^{-1})$.
Combined with \eqref{eqn:41}, this yields the desired result.
\qed

\medskip

We record here a useful bound that was used in the previous proof:
Recall from Section \ref{sec:notation} that $C_{\vt}$ denotes the contour function of $\vt\in\bM_{n, 0}$.
Then
\begin{align}\label{eqn:47}
&
\#\big\{(i_1,\ldots, i_{2s})\ :\ i_{2j}\in\fB(C_{\vt}; i_{2j-1})\text{ for } 1\leq j\leq s\, ,\ 
i_{\ell}\text{-s NAD}  \big\}\notag\\
&\hskip40pt\leq
s!\times
\#\big\{\big(i_1,\ldots, i_{2s}, k_1,\ldots, k_{2s}\big)\in\bfac(\vt, s)\ :\  i_{\ell}\text{-s NAD}
\big\}\notag\\
&\hskip80pt\leq
c n^s\cdot\big(\max_{y}L\big(C_{\vt};\, 2n, y\big)\big)^{s-1}=:\alpha_n(\vt)\, ,
\end{align}
where $c$ is a constant that depends only on $s$. 
We are now ready for the

\medskip

\noindent{\bf Proof of Lemma \ref{lem:2}:}
Recall the definition of $(T_n^{\bbf},\, \Xi_n^{\bbf})$ from Proposition \ref{prop:bf-scaling-limit}.
Clearly,
\[
\pr\big(T_n^{\bbf}=\vt\big)
=
\frac{\#\bfac(\vt, s)}{\Sigma\bfac}\, ,\ \ \ \ \vt\in\bM_{n, 0}\, ,
\]
where $\Sigma\bfac=\sum_{\vt'\in\bM_{n, 0}}\#\bfac(\vt', s)$.
Then
\begin{align*}
&
\sum_{\vt\in\bM_{n, 0}} \big|\pr\big(T_n^{\circ}=\vt\big)-\pr\big(T_n^{\bbf}=\vt\big)\big|
=
\sum_{\vt\in\bM_{n, 0}} \bigg|\frac{B(C_{\vt})^s}{\sum_{C\in\fC_n}B(C)^s}
-
\frac{\#\bfac(\vt, s)\times s!}{\Sigma\bfac\times s!}\bigg|
\notag\\
&
\leq
\sum_{\vt\in\bM_{n, 0}} \bigg|\frac{B(C_{\vt})^s-\#\bfac(\vt, s)\times s!}{\sum_{C\in\fC_n}B(C)^s}\bigg|
+
\sum_{\vt\in\bM_{n, 0}}
\#\bfac(\vt, s)\times s!\cdot\bigg|\frac{1}{\sum_{C\in\fC_n}B(C)^s}-\frac{1}{\Sigma\bfac\times s!}\bigg|
\notag\\
&\hskip40pt
\leq
2\sum_{\vt\in\bM_{n, 0}} \bigg|\frac{B(C_{\vt})^s-\#\bfac(\vt, s)\times s!}{\sum_{C\in\fC_n}B(C)^s}\bigg|
=
2\cdot
\frac{\bE\big|B(C_n)^s-\#\bfac(T_n, s)\times s!\big|}{\bE\big[B(C_n)^s\big]}
\to 0 \, ,
\end{align*}
where the last step uses Lemma \ref{lem:3}, \eqref{eqn:29}, and \eqref{eqn:44}.
Let us assume that $T_n^{\bbf}$ and $T_n^{\circ}$ are coupled in a way so that
\begin{align}\label{eqn:45}
\pr\big(T_n^{\bbf}\neq T_n^{\circ}\big)\to 0\, ,\ \ \ \text{ as }\ \ \ n\to\infty\, .
\end{align}

Let 
\[
\Xi_n^{\bbf}=\big(I_{n, j}^{\bbf}, K_{n, j}^{\bbf}\, ;\, 1\leq j\leq 2s\big)\, ,\ \ 
\pr_{\vt}^{\bbf}\big(\cdot\big)=\pr\big(\cdot\ |\ T_n^{\bbf}=\vt\big)\, ,\ \text{ and }\
\pr_{\vt}^{\circ}\big(\cdot\big)=\pr\big(\cdot\ |\ T_n^{\circ}=\vt\big)\, .
\]
Using \eqref{eqn:47}, we see that
\begin{gather}\label{eqn:48}
\pr_{\vt}^{\bbf}\big(I_{n, j}^{\bbf}\text{-s NAD}\big)\leq\frac{\alpha_n(\vt)}{\#\bfac(\vt, s)}=:\beta_n^{(1)}(\vt)\, , \text{ and }\
\pr_{\vt}^{\circ}\big(I_{n, j}^{\circ}\text{-s NAD}\big)\leq\frac{\alpha_n(\vt)}{B(C_{\vt})^s}=:\beta_n^{(2)}(\vt)\, .
\end{gather}
Let 
\[
\Gamma_{n, s}(\vt)=\big\{(i_1,\ldots, i_{2s})\in[2n]^{2s}\, :\, 
\exists\ k_1,\ldots, k_{2s}\text{ such that }
(i_j, k_j\, ;\, 1\leq j\leq 2s)\in\bfac(\vt, s)
\big\}\, .
\]
Write $\mvI_n^{\bbf}=(I_{n, 1}^{\bbf},\ldots,I_{n, 2s}^{\bbf})$ and similarly define $\mvI_n^{\circ}$.
Then for any $\mvj=(j_1,\ldots,j_{2s})\in\Gamma_{n,s}(\vt)$ with $j_{\ell}$-s all distinct,
\begin{align*}
\bigg|
\pr_{\vt}^{\circ}\big(\sort(\mvI_n^{\circ})=\mvj\big)
-
\pr_{\vt}^{\bbf}\big(\mvI_n^{\bbf}=\mvj\big)
\bigg|
=\frac{s!}{B(C_{\vt})^s}-\frac{1}{\#\bfac(\vt, s)}\, ,
\end{align*}
where
$\sort(\mvI_n^{\circ})
=
\big(I_{n,\pi(1)}^{\circ}, I_{n,\pi(1)+1}^{\circ}, I_{n,\pi(3)}^{\circ}, I_{n,\pi(3)+1}^{\circ},\ldots, I_{n,\pi(2s-1)}^{\circ}, I_{n,\pi(2s-1)+1}^{\circ}\big)$
for a permutation $\pi$ of $\{1, 3, \ldots, 2s-1\}$ such that $\sort(\mvI_n^{\circ})\in\Gamma_{n,s}(\vt)$.
Hence,
\begin{align}\label{eqn:49}
\sum_{\substack{\mvj\in\Gamma_{n, s}(\vt)\\ j_{\ell}\text{-s all distinct}}}
\bigg|
\pr_{\vt}^{\circ}\big(\sort(\mvI_n^{\circ})=\mvj\big)
-
\pr_{\vt}^{\bbf}\big(\mvI_n^{\bbf}=\mvj\big)
\bigg|
\leq
\frac{s!\times\#\bfac(\vt, s)-B(C_{\vt})^s}{B(C_{\vt})^s}
=:
\beta_n^{(3)}(\vt)
\, .
\end{align}
Using \eqref{eqn:48} and \eqref{eqn:49}, we see that we can sample $\mvI_n^{\bbf}$ and $\mvI_n^{\circ}$ in a way so that
\[
\pr\big(\sort(\mvI_n^{\circ})\neq\mvI_n^{\bbf}\ \big|\ T_n^{\circ}=T_n^{\bbf}=\vt\big)
\leq
\sum_{j=1}^3\beta_n^{(j)}(\vt)=:\beta_n(\vt)\, .
\]
Combined with \eqref{eqn:45}, we get a coupling of $G_{n,s}^{\circ}$ and $G_{n,s}^{\bbf}$ such that for any $\eps>0$,
\begin{align}\label{eqn:50}
\pr\big(G_{n,s}^{\circ}\neq G_{n,s}^{\bbf}\big)
&\leq
\pr\big(T_n^{\circ}\neq T_n^{\bbf}\big)
+
\sum_{\vt : \beta_n(\vt)<\eps}\beta_n(\vt)\cdot\pr\big(T_n^{\circ}=T_n^{\bbf}=\vt\big)
+
\pr\big(\beta_n(T_n^{\circ})\geq\eps\big)
\notag\\
&\leq
o(1)
+
\eps
+
\pr\big(\beta_n(T_n^{\circ})\geq\eps\big)\, .
\end{align}
To complete the proof, it is enough to show that $\pr\big(\beta_n(T_n^{\circ})\geq\eps\big)\to 0$.
Note that
\begin{align}\label{eqn:51}
\pr\big(\beta_n(T_n^{\circ})\geq\eps\big)
=
\frac{\bE\big[\ind\{\beta_n(T_n)\geq\eps\}\cdot B(C_n)^s\big]}{\bE\big[B(C_n)^s\big]}
\leq 
\frac{\pr\big(\beta_n(T_n)\geq\eps\big)^{1/2}\cdot\bE\big[ B(C_n)^{2s}\big]^{1/2}}{\bE\big[B(C_n)^s\big]}\, .
\end{align}
In view of \eqref{eqn:29} and \eqref{eqn:44}, it is enough to show that $\pr\big(\beta_n(T_n)\geq\eps\big)\to 0$, which is a simple consequence of \eqref{eqn:25}, Lemma \ref{lem:3}, and \eqref{eqn:29}.
\qed

\medskip

\subsection{Proof of Proposition \ref{prop:df-scaling-limit}}\label{sec:proof-df-scaling-limit}
We tailor the argument in the proof of Proposition \ref{prop:bf-scaling-limit} to the depth-first setting.
We will work with a fixed $s\geq 1$ throughout this section.

\begin{figure}
	\centering
	\includegraphics[trim=.85cm 20cm 3.5cm 3cm, clip=true, angle=0, scale=.8]{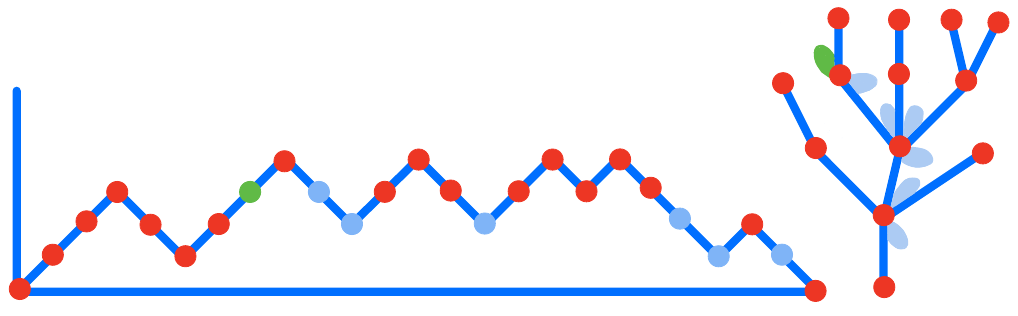}
	\captionof{figure}{Left: $f$ with the point on the graph of $f$ corresponding to $i=7$ colored green, and points that correspond to $\fD(f; 7)\setminus\{7\}$ colored light blue. 
	Right: the corresponding plane tree with the relevant corners colored.}
	\label{fig:8}
\end{figure}

For $a<b$  and a continuous excursion $e$ on $[a, b]$, $t\in[a, b]$, and $y\in [0, e(t)]$, define
\[
\fR\big(e; t, y\big):=\big\{u\in[t, b]\ :\ e(u)=y\, ,\ \nexists v\in[t, u]\text{ such that }e(v)<y\big\}\, .
\]
Let $\fC_{n}$ be as in \eqref{eqn:52}.
For $f\in\fC_{n}$ and $i=0, 1,\ldots, 2n$, let
\begin{gather*}
\fD(f; i):= \bigcup_{j=1}^{f(i)}\fR\big(f; i, j\big)\, ,\
D(f; i):=\#\, \fD(f; i)\, , \ \text{ and }\ D(f)=\sum_{i=0}^{2n}D(f, i)\, .
\end{gather*}
An illustration of $\fD(f; i)$ is given in Figure \ref{fig:8}.
In terms of the plane tree whose contour function is $f$, for $1\leq i\leq 2n-1$, 
$\fD(f; i)\setminus\{i\}$ corresponds to the set of corners that (i) are incident to vertices on the ancestral line of the vertex to which the $i$-th corner is incident, (ii) appear strictly after the $i$-th corner in the contour order, and (iii) are distinct from the corner associated with the root edge.
Note the connection with the set of depth-first admissible corners $\dfac(\cdot, \cdot)$.

Let $C_n^{\bullet}$ be distributed as
\begin{align}\label{eqn:53}
\pr\big(C_{n}^{\bullet}=f\big)
=
\frac{\big(D(f)\big)^s}{\sum_{\phi\in\fC_{n}}\big(D(\phi)\big)^s}\, , \ \ \ \  f\in\fC_{n}\, .
\end{align}
Conditional on $C_n^{\bullet}$, sample i.i.d. random variables 
$\big(I_{n, 2j-1}^{\bullet}\, ,\, I_{n,2j}^{\bullet}\big)$, $1\leq j\leq s$, where
\begin{align}\label{eqn:54}
\pr\big(I_{n, 1}^{\bullet}=i\big)
=
D(C_n^{\bullet}; i)/D(C_n^{\bullet})\, , \ \ i=0, 1,\ldots, 2n\, ,
\end{align}
and conditional on $I_{n, 1}^{\bullet}$, $I_{n, 2}^{\bullet}$ is uniform over $\fD(C_n^{\bullet}; I_{n, 1}^{\bullet})$.

Let $T_n^{\bullet}$ be the plane tree whose contour function is $C_n^{\bullet}$.
Let $G_{n, s}^{\bullet}$ be the rooted metric measure space obtained by adding an edge in $T_n^{\bullet}$ between the $I_{n, 2j-1}^{\bullet}$-th corner and the $I_{n,2j}^{\bullet}$-th corner of $T_n^{\bullet}$, $1\leq j\leq s$, endowing the resulting space by the graph distance, rooting it at the root vertex of $T_n^{\bullet}$, and assigning probability $1/n$ to every vertex of $T_n^{\bullet}$ except the root vertex. 
The next two lemmas combined complete the proof of Proposition \ref{prop:df-scaling-limit}.

\begin{lem}\label{lem:4}
	We have, 
	\[
	\frac{1}{\sqrt{n}}\cdot G_{n, s}^{\bullet}\weakc \frac{1}{\sqrt{2}}\cdot\cH_{(s)},\ \ \text{ as }\ \ n\to\infty
	\]
	w.r.t. the pointed GHP topology.
\end{lem}
\begin{lem}\label{lem:5}
	There exists a coupling of $G_{n, s}^{\bullet}$ and $G_{n, s}^{\df}$ such that
	\[
	\pr\big(G_{n, s}^{\bullet}\neq G_{n, s}^{\df}\big)\to  0,\ \ \text{ as }\ \ n\to\infty\, .
	\]
\end{lem}

We need some control over the functional $D(C_n^{\bullet}\, ;\, \cdot)$ in order to prove Lemma \ref{lem:4}.
To this end, for $\vt\in\bM_{n, 0}$, let $\big(S_{\vt}(\ell)\, ,\, 0\leq \ell\leq n+1\big)$ be the {\L}ukasiewicz path (see, e.g., \cite{legall-survey} for definition) of $\vt$.
For $i\in\{0, 1,\ldots, 2n-1\}$,
consider the vertex of $\vt$ to which the $i$-th corner of $\vt$ is incident, and 
suppose that the {\L}ukasiewicz path of $\vt$ explores this vertex at the $\ell_i$-th step.
Further, denote by $\zeta_{\vt}(i)$ the degree of this vertex in $\vt$.
Then
\begin{align}\label{eqn:55}
S_{\vt}(\ell_i)-\zeta_{\vt}(i)
\leq
D(C_{\vt}; i)-C_{\vt}(i)
\leq
S_{\vt}(\ell_i)+1\, .
\end{align}
As in the previous section, let $C_n$ be uniformly distributed over $\fC_n$, and let $T_n$ be the plane tree whose contour function is $C_n$.
Recall that the tree obtained from $T_n$ by deleting its root edge has the same distribution as a uniform plane tree on $n$ vertices.
This observation together with \cite[Theorem 3]{marckert-mokkadem} for the case of $\mathrm{Geometric}(1/2)$ offspring distribution gives
\[
n^{-1/2}\max_{i}\big|S_{T_n}(\ell_i)-C_n(i)\big|\weakc 0\, .
\]
Combined with \eqref{eqn:55} and the fact $\max_i\zeta_{T_n}(i)=O_P(\log n)$, this yields
\begin{align}\label{eqn:56}
n^{-1/2}\max_{i}\big|D(C_n; i)-2C_n(i)\big|\weakc 0\, .
\end{align}
Further, it follows from 
\cite[Equation 32 and Theorem 1.2]{addarioberry-devroye-janson} that\footnote{
\cite[Equation 32]{addarioberry-devroye-janson} is stated for the maximum of the breadth-first queue, whereas we are dealing with the maximum of the depth-first queue.
However, it is well-known that these two have the same distribution.
} 
for all $x>0$,
\[
\pr\big(\max_{i}|S_{T_n}(i)|\geq x\sqrt{n}\big)\leq ce^{-c' x^2}\, ,\ \text{ and }\
\pr\big(\max_{i}C_n(i)\geq x\sqrt{n}\big)\leq c e^{-c' x^2} \, ,
\]
which combined with \eqref{eqn:55} yields
\begin{align}\label{eqn:57}
\pr\big(\max_{i} D\big(C_n; i\big)\geq x\sqrt{n}\big)\leq ce^{-c' x^2}\, ,\ \ \text{ for all }\ \ x>0\, .
\end{align}

\medskip

\noindent{\bf Proof of Lemma \ref{lem:4}:} 
First, using \eqref{eqn:56} and \eqref{eqn:24}, we see that
\begin{align}\label{eqn:58}
\frac{D(C_n)}{(2n)^{3/2}}\weakc 2\int_{0}^{1}\ve(t) dt\, ,\ \ \text{ as }\ \ n\to\infty\, ,
\end{align}
jointly with the convergence $\bar C_n\weakc\ve$.
Also, it follows from \eqref{eqn:57} and the bound $D(C_n)\leq 2n\max_i D(C_n; i)$ that for any $s\geq 1$,  the sequence of random variables 
\begin{align}\label{eqn:60}
\big\{(2n)^{-3s/2}D(C_n)^s\big\}_{n\geq 1}\ \text{ is uniformly integrable.}
\end{align}
Let 
\[
\bar C_n^{\bullet}(t):=(2n)^{-1/2}C_n^{\bullet}(2nt)\, ,\ \ 0\leq t\leq 1\, .
\]
Then for any bounded continuous $\phi:C[0, 1]\to\bR$,
\begin{align*}
\bE\big[\phi\big(\bar C_n^{\bullet}\big)\big]
=
\frac{\bE\big[\phi\big(\bar C_n\big)D(C_n)^s\big]}{\bE\big[D(C_n)^s\big]}
\ \tonn\
\frac{\bE\big[\phi(\ve)\big(\int_0^1\ve(t)dt\big)^s\big]}{\bE\big[\big(\int_0^1\ve(t)dt\big)^s\big]}
=
\bE\big[\phi(\esdf)\big]\, ,
\end{align*}
where the second step uses \eqref{eqn:24}, \eqref{eqn:58}, and \eqref{eqn:60}.
Repeating the argument, we can assume that we are working in a space where 
\begin{gather}\label{eqn:59}
\bar C_n^{\bullet}\convas\esdf\, ,\ 
\frac{1}{\sqrt{n}}\max_{i}\big|D(C_n^{\bullet}; i)-2C_n^{\bullet}(i)\big|\convas 0\, ,\  \ \text{ and }\ \
\frac{D(C_n^{\bullet})}{(2n)^{3/2}}\convas  2\int_{0}^{1}\esdf(t)dt\, .
\end{gather}
Using \eqref{eqn:54}, we can further assume that on this space,
\begin{gather}\label{eqn:61}
(2n)^{-1}I_{n,1}^{\bullet}\convas u_1\, ,
\end{gather}
where $u_1$ is as in Construction \ref{constr:H-s}.

Sample $I_{n, 2}^{\bullet}$ in the following two steps:
Conditional on $C_n^{\bullet}, I_{n, 1}^{\bullet}$, first sample $J_n^{\bullet}$ according to \greg{ the } law
\begin{align}\label{eqn:62}
\pr\big(J_n^{\bullet}=j\ \big|\ C_n^{\bullet}, I_{n, 1}^{\bullet}\big)
=
\frac{\#\fR\big(C_n^{\bullet}\, ;\,  I_{n, 1}^{\bullet}, j\big)}{D\big(C_n^{\bullet} , I_{n, 1}^{\bullet}\big)}\, ,\ \ \
j=1,\ldots, C_n^{\bullet}(I_{n,1}^{\bullet})\, ,
\end{align}
and then let 
$I_{n,2}^{\bullet}\sim
\mathrm{Uniform}\big(\fR\big(C_n^{\bullet}\, ;\,  I_{n, 1}^{\bullet}, J_n^{\bullet}\big)\big)$.
Let 
$i_j=\min\, \fR\big(C_n^{\bullet}\, ;\,  I_{n, 1}^{\bullet}, j\big)$ 
for $j=1,\ldots, C_n^{\bullet}(I_{n,1}^{\bullet})$.
Then 
\[
C_n^{\bullet}(i_j)=j\, , \ \ \text{ and }\ \
D(C_n^{\bullet} ;\, i_j)-D(C_n^{\bullet} ;\, i_{j-1})=\#\fR\big(C_n^{\bullet}\, ;\,  i_j, j\big)\, .
\]
Consequently, using \eqref{eqn:62} and the second convergence in \eqref{eqn:61}, we see that
\[
\pr\big(J_n^{\bullet}\leq j\ \big|\ C_n^{\bullet}, I_{n, 1}^{\bullet}\big)
=
\frac{D(C_n^{\bullet} ;\, i_j)}{D\big(C_n^{\bullet} , I_{n, 1}^{\bullet}\big)}
=
\frac{2j+o(\sqrt{n})}{2C_n^{\bullet}(I_{n, 1}^{\bullet})+o(\sqrt{n})}\, ,\ \ \
j=1,\ldots, C_n^{\bullet}(I_{n,1}^{\bullet})\, ,
\]
where the $o(\sqrt{n})$ term is uniform over $j$.
Thus, using \eqref{eqn:61} and the first convergence in \eqref{eqn:59}, we can assume that 
\begin{align}\label{eqn:63}
(2n)^{-1/2}\cdot J_n^{\bullet}\convas z_1
\end{align}
together with \eqref{eqn:59} and \eqref{eqn:61}, where $z_1$ is as in Construction \ref{constr:H-s}.
Let $v_1$ be as in Construction \ref{constr:H-s} as well.
Then the set $\fR\big(\esdf\, ;\, u_1, z_1\big)$ is the singleton $\{v_1\}$ almost surely, and hence the first convergence in \eqref{eqn:59}, \eqref{eqn:61}, and \eqref{eqn:63} imply that
\[
(2n)^{-1}\min\, \fR\big(C_n^{\bullet}\, ;\,  I_{n, 1}^{\bullet}, J_n^{\bullet}\big)\convas v_1\, ,\ \ 
\text{ and }\ \ 
(2n)^{-1}\max\, \fR\big(C_n^{\bullet}\, ;\,  I_{n, 1}^{\bullet}, J_n^{\bullet}\big)\convas v_1\, .
\]
In particular, $(2n)^{-1}I_{n, 2}^{\bullet}\convas v_1$.
Repeating the same argument, we can assume that
\[
(2n)^{-1}\big(I_{n, 2j-1}^{\bullet}\, ,\, I_{n, 2j}^{\bullet}\big)
\convas
\big(u_j, v_j\big)\, ,\ \ 1\leq j\leq s\, .
\]
Now the proof can be completed by following the argument given after \eqref{eqn:34}.\qed

\medskip

\noindent{\bf Proof of Lemma \ref{lem:5}:}
The proof follows the same steps as in the breadth-first case.
First note that we can use the argument in the proof of Lemma \ref{lem:3} to prove the following analogous result:
\[
n^{-3s/2}\cdot\bE\big|s!\times\#\dfac(T_n, s)- D(C_n)^s\big|\to 0\ \ \text{ as }\ \ n\to\infty\, ;
\]
here the bound in \eqref{eqn:57} plays the role analogous to that of \eqref{eqn:25} in the breadth-first setting.
Using the above convergence, we can construct, similar to \eqref{eqn:45}, 
a coupling of $T_n^{\df}$ and $T_n^{\bullet}$ such that
$\pr\big(T_n^{\df}\neq T_n^{\bullet}\big)\to 0$ as $n\to\infty$.
Then the rest of the argument in the proof of Lemma \ref{lem:2} can be mimicked to complete the proof.
We omit the details to avoid repetition. 
\qed

\subsection{Proof of Corollary \ref{cor:radius-Hs}}
\greg{By Jeulin's local time identity \cite{jeulin,AMP},
\begin{align}\label{eqn:99}
\big(\eta(\ve; 1,y/2),\, y\geq 0\big)
\equald
\big(2\cdot\ve(\tau^{-1}(y)),\, y\geq 0\big)\, ,\
\text{ where }\
\tau(t)=\int_0^t\frac{d u}{\ve(u)}\, \, ,\quad t\geq 0\, ,
\end{align}
and $\tau^{-1}(y)=\sup\big\{t\in [0, 1]\, :\, \tau(t)\leq y\big\}$ for $y\geq  0$. 
Now,
\begin{align}\label{eqn:101}
\int_0^{\infty}\eta(\ve; 1, y)^2 dy
=\frac{1}{2}\int_0^{\infty}\eta(\ve; 1, y/2)^2 dy
\equald 2\int_0^{\infty}\big(\ve(\tau^{-1}(y))\big)^2 dy
=2\int_0^1\ve(t)dt\, ,
\end{align}
where the second step follows from \eqref{eqn:99}, and the last step follows if we use the substitution $y=\tau(t)$.
Further, since $2\cdot\|\ve\|_\infty=\inf\big\{y>0:\eta(\ve; 1,y/2)=0\big\}$ almost surely, \eqref{eqn:99} implies 
\begin{align}\label{eqn:100}
2\|\ve\|_\infty
\equald
\tau(1)\, 
\end{align}
jointly with the equality in distribution in \eqref{eqn:101}.
Consequently, for every bounded continuous $\phi:\bR\to\bR$, 
\begin{align*}
\bE\big[\phi\big(\radi(\cH_{(s)})\big)\big]
&
=\bE\big[\phi\big(2\|\esbf\|_{\infty}\big)\big]
=\frac{\bE\big[\phi\big(2\|\ve\|_{\infty}\big)\big(\int_0^\infty \eta(\ve; 1,y)^2 dy\big)^s\big]}{\bE\big[\big(\int_0^\infty \eta(\ve; 1,y)^2 dy\big)^s\big]}\\
&
=\frac{\bE\big[\phi\big(\int_0^1 dt/\ve(t)\big)\big(\int_0^1 \ve(t) dt\big)^s\big]}{\bE\big[\big(\int_0^1 \ve(t)dt\big)^s\big]}
=\bE\big[\phi\big(\int_0^1 \frac{dt}{\esdf(t)}\big)\big]
\, ,
\end{align*}
where the first step uses \chng{Corollary~\ref{cor:bf=df-Hs}~(i)}, and the third step uses \eqref{eqn:101} and \eqref{eqn:100}.
This completes the proof.
}

\subsection{Proof of Theorem \ref{thm:H-n-s-width-convergence}}\label{sec:H-n-s-width-convergence-proof}
Let us first state the result concerning convergence of the height profile of $H_{n, 0}$.
\begin{thm}\label{thm:total-local-time}
Let $Z_{n, 0}(\ell)$ denote the number of vertices in $H_{n, 0}$ at distance $\ell$ from the root, $\ell=0,1,\ldots$.
Let $\bar Z_{n, 0}(r)=n^{-1/2}\cdot Z_{n,0}\big(\lfloor r\sqrt{n}\rfloor\big)$, $r\geq 0$.
Then 
\begin{align}\label{eqn:72}
\Big(\bar Z_{n, 0}(r),\ r\geq 0\Big)
\weakc
\Big(\eta\big(2\ve; 1, r\big),\ r\geq 0\Big)
=
\Big(\frac{1}{2}\eta\big(\ve; 1, r/2\big),\ r\geq 0\Big)\, ,\ \ \text{ as }\ \ n\to\infty\, ,
\end{align}
w.r.t. Skorohod $J_1$ topology on $\bD([0,\infty)\ :\ \bR)$.
\end{thm}
As mentioned before, Theorem \ref{thm:total-local-time} is a special case of \cite[Theorem 1.1]{drmota-gittenberger}. 
We will make use of this result in our proof.

For $s\geq 1$, we can explore any $G\in\bH_{n, s}$ in a breadth-first manner and get a breadth-first spanning tree as follows: 
Add an extra vertex labeled `$0$' and connect it to the root of $G$ via an edge, and declare this edge oriented away from $0$ to be the root edge.
Use the vertex labels to endow the neighbors of every vertex in the resulting graph with a circular order (where our convention is to go from the smallest label to the highest). 
Then we can explore this map in a breadth-first fashion as in Section \ref{sec:bf-df-search}.
From the resulting tree we delete the vertex $0$ and the edge incident to it, and root this tree at the root vertex of $G$.
Thus, we end up with a tree $\vt\in \bH_{n, 0}$.

Now there exist $(i_k , j_k)\in[n]^2$, $1\leq k\leq s$, such that 
\begin{align}\label{eqn:66}
0\leq \hght(i_k; \vt)-\hght(j_k; \vt)\leq 1\, ,
\end{align}
and
$G$ can be recovered by adding an edge in $\vt$ between the vertices $i_k$ and $j_k$ for each $1\leq k\leq s$.
Think of $\vt$ as a plane tree by using the vertex labels and consider the plane embedding.
Then in this embedding, for each $k$, 
the vertex $j_k$ will appear after the vertex $i_k$ in the \greg{depth-first} order.
%\begin{align}\label{eqn:67}
%\text{the vertex }\ j_k\ \text{ will appear after the vertex }\ i_k\ \text{ in the contour order.}
%\end{align}
(These conditions are similar to \eqref{eqn:16} and the condition `$i_{2j-1}\leq i_{2j}$' appearing right above it.)

Thus, there is an asymmetry in the roles of $i_k$ and $j_k$.
This does not cause a problem in the proof of Proposition \ref{prop:bf-scaling-limit} because of the convergence of the second coordinate in \eqref{eqn:24}, i.e., the convergence,  {\it as a function of two variables}, of the discrete local time of the contour function to a continuous limit.
The analogue of this result for uniform labeled trees is not available in the literature, which poses a problem in working with $\vt$ directly.
This will be discussed further in Section \ref{sec:discussion}.
We, however, do have convergence of the total local time (Theorem \ref{thm:total-local-time}).
To go around the above issue by making use of Theorem \ref{thm:total-local-time}, we define a new tree $\bar\vt$ by applying a kind of symmetrization to $\vt$.

If there exist $1\leq k\neq q\leq s$ such that \greg{$\big|\hght(i_k;\vt)-\hght(i_q; \vt) \big|\leq 1$}, set $\bar\vt$ to be the empty tree $\emptyset$.
Otherwise, do the following for $1\leq k\leq s$:
If $\hght(i_k; \vt)=\hght(j_k; \vt)$, do nothing,
and if $\hght(i_k; \vt)=\hght(j_k; \vt)+1$,
then with probability $1/2$ do nothing, and with probability $1/2$
add an edge in $\vt$ between $i_k$ and $j_k$ and delete the edge between $i_k$ and its parent in $\vt$; call this operation a `swap.'
Denote the resulting tree by $\bar\vt$.
We set 
\[
\sbf(G)=\vt\ \ \ \text{ and }\ \ \ \overline\sbf(G)=\bar\vt\; .
\]
Note that the swapping operations above commute, i.e., it does not matter in what order they are done.
Note also that for any $k\geq 1$, the number of vertices in $G$ that are at distance $k$ from the root is same as the number of vertices at height $k$ in $\sbf(G)$, and hence the same is true of $\overline\sbf(G)$ whenever $\overline\sbf(G)\neq\emptyset$.

To fix ideas, first consider $s=1$.
In this case, $\overline\sbf(H_{n, 1})$ cannot be the empty tree $\emptyset$.
For any $\bar\vt\in\bH_{n, 0}$, the event $\big\{\overline\sbf(H_{n, 1})=\bar\vt\big\}$ holds true iff there exist $i, j\in[n]$ such that 
$i\neq j$ and $j$ is not the parent of $i$ in $\bar\vt$, 
$H_{n, 1}$ is the graph obtained by adding an edge in $\bar\vt$ between $i$ and $j$, and exactly one of the following happens:

\vskip5pt

\begin{inparaenumi}
\noindent\item 
$\hght(j; \bar\vt)=\hght(i; \bar\vt)$ and $i<j$.
(In this case, $\sbf(H_{n, 1})=\overline\sbf(H_{n, 1})=\bar\vt$.)

\vskip5pt

\noindent\item 
Think of $\bar\vt$ as a plane tree %by using the vertex labels, 
\greg{in which the children of a vertex are ordered from left to right in order of increasing label.} 
%and consider its plane embedding.
Then $\hght(j; \bar\vt)=\hght(i; \bar\vt)-1$, and
$j$ appears after $i$ in the contour order of $\bar\vt$,
and swapping did not take place in going from $\sbf(H_{n, 1})$ to $\overline\sbf(H_{n, 1})$.
(In this case also, $\sbf(H_{n, 1})=\overline\sbf(H_{n, 1})=\bar\vt$.)

\vskip5pt

\noindent\item
$\hght(j; \bar\vt)=\hght(i; \bar\vt)-1$, 
$j$ appears before $i$ in the contour order of $\bar\vt$,
and swapping took place in going from $\sbf(H_{n, 1})$ to $\overline\sbf(H_{n, 1})$.
(In this case, $\sbf(H_{n, 1})$ is the tree obtained from $\bar\vt$ by adding an edge between $i$ and $j$ and deleting the edge between $i$ and its parent in $\bar\vt$.)
\end{inparaenumi}

\vskip5pt

Then we see from the above discussion that for $\bar\vt\in\bH_{n, 0}$,
\begin{align}\label{eqn:65}
\pr\big(\overline\sbf(H_{n, 1})=\bar\vt\big)
=
\frac{1}{\#\bH_{n, 1}}\cdot
\sum_{\ell\geq 1}\Bigg[
\Bigg(
\begin{array}{c}
\hskip-4pt z(\bar\vt; \ell)\hskip-7pt\phantom{.} \\
\hskip-4pt 2\hskip-7pt \phantom{.}
\end{array}
\Bigg)
+
\frac{1}{2}\cdot z(\bar\vt; \ell)\cdot\Big(z\big(\bar\vt; \ell-1\big)-1\Big)\Bigg]
=:\frac{W_1(\bar\vt)}{\#\bH_{n, 1}}
\, ,
\end{align}
where $z(\cdot ; \cdot)$ is as defined in Section \ref{sec:notation}.
As observed before, $H_{n, 1}$ and $\overline\sbf(H_{n, 1})$ have the same distance profile.
Thus, using \eqref{eqn:65}, we see that for any bounded continuous $\phi:\bD[0, \infty)\to\bR$,
\begin{align}\label{eqn:70}
\bE\big[\phi\big(\bar Z_{n, 1}\big)\big]
=
\frac{\#\bH_{n, 0}}{\#\bH_{n, 1}}\cdot
\bE\big[\phi\big(\bar Z_{n, 0}\big)W_1(H_{n, 0})\big]\, .
\end{align}
From Theorem \ref{thm:total-local-time}, it follows that
\begin{align}\label{eqn:71}
n^{-3/2}\cdot W_1(H_{n, 0})\weakc
\frac{1}{2}\int_{0}^{\infty}\eta\big(\ve; 1, y\big)^2 dy
\end{align}
jointly with \eqref{eqn:72}.
Further, for any $s\geq 0$,
\begin{align}\label{eqn:73}
n^{1-n-3s/2}\cdot(\#\bH_{n, s})
\tonn
\frac{1}{s!}\cdot\E\bigg[\bigg(\int_0^1\ve(x)dx\bigg)^s\bigg]
=
\frac{1}{2^s s!}\cdot\E\bigg[\bigg(\int_{0}^{\infty}\eta\big(\ve; 1, y\big)^2 dy\bigg)^s\bigg]\, ,
\end{align}
where the first step follows from the results of \cite{wright1977number, spencer-count}, and the second step uses \eqref{eqn:74}.
Finally, using the bound
$W_1(\bar\vt)\leq n\cdot\max_{\ell}z(\bar\vt, \ell)$
and applying Theorem \ref{thm:conditioned-gw-tree} with $\mathrm{Poisson}(1)$ offspring distribution, we see that 
$n^{-3/2}\cdot W_1(H_{n, 0})$ is a sequence of uniformly integrable random variables.
Combining this last observation with \eqref{eqn:70}, \eqref{eqn:71}, and \eqref{eqn:73}, and Theorem \ref{thm:total-local-time}, we get
\begin{align}\label{eqn:75}
\bE\big[\phi\big(\bar Z_{n, 1}\big)\big]
\tonn
\frac{\bE\big[\phi\big(\eta\big(2\ve; 1, \cdot\big)\big) \int_{0}^{\infty}\eta\big(\ve; 1, y\big)^2 dy\big]}{\bE\big[\int_{0}^{\infty}\eta\big(\ve; 1, y\big)^2 dy\big]}
=\bE\big[\phi\big(\eta\big(2\ve_{(1)}^{\mathrm{BF}}; 1, \cdot\big)\big)\big]\, .
\end{align}
This completes the proof for $s=1$.

For $s\geq 2$, we need some control on $\pr\big(\overline\sbf(H_{n, s})=\emptyset\big)$.
To this end, note that for any $\vt\in\bH_{n, 0}$, the number of tuples $(i_k, j_k)$, $1\leq k\leq s$, that satisfy \eqref{eqn:66}, and 
\greg{$\big|\hght(i_k;\vt)-\hght(i_q; \vt) \big|\leq 1$}
for some $1\leq k\neq q\leq s$ is upper bounded by
\begin{align}\label{eqn:68}
\gamma(\vt):=c n^{s-1}\big[\max_\ell z(\vt; \ell)\big]^{s+1}\, ,
\end{align}
where $c>0$ depends only on $s$\greg{, for instance we can choose $c=3\cdot s(s-1)\cdot 2^s$}. 
Consequently,
\begin{align}\label{eqn:69}
\pr\big(\overline\sbf(H_{n, s})=\emptyset\big)
&\leq
\frac{1}{\#\bH_{n, s}}\cdot\sum_{\vt\in\bH_{n, 0}}\gamma(\vt)\\
&=cn^{s-1}\cdot(\sqrt{n})^{s+1}\cdot\bigg(\frac{\#\bH_{n, 0}}{\#\bH_{n, s}}\bigg)\cdot\bE\bigg[\left(\frac{\max_\ell z(H_{n, 0}; \ell)}{\sqrt{n}}\right)^{s+1}\bigg]
\leq\frac{c'}{\sqrt{n}}\, ,\notag
\end{align}
where the last step follows if we apply Theorem \ref{thm:conditioned-gw-tree} with $\mathrm{Poisson}(1)$ offspring distribution, and use \eqref{eqn:73}.

Now for any $\bar\vt\in\bH_{n, 0}$, the event $\big\{\overline\sbf(H_{n, s})=\bar\vt\big\}$ holds true iff 
there exist $i_k, j_k\in[n]$, $1\leq k\leq s$, such that the following hold:
\begin{inparaenumi}
\item
either 
$\hght(i_k; \bar\vt)=\hght(j_k; \bar\vt)$ and $i_k<j_k$ or
$\hght(i_k; \bar\vt)=\hght(j_k; \bar\vt)+1$ and $j_k$ is not the parent of $i_k$ in $\bar\vt$,
\item
\greg{$\hght(i_{k+1}; \bar\vt)\geq\hght(i_k; \bar\vt)+2$} for  $k=1,\ldots, s-1$,
\item
$H_{n, s}$ is the graph obtained from $\bar\vt$ by placing an edge between $i_k$ and $j_k$ for $1\leq k\leq s$, and
\item
the unique swaps (if any) needed to go from $\sbf(H_{n,s})$ to $\overline\sbf(H_{n,s})$ were performed.
\end{inparaenumi}

\vskip5pt

\noindent Writing $\sum_1$ for sum over all $\ell_1, \ldots, \ell_s$ such that \greg{$\ell_{k+1}\geq\ell_k +2$}, $1\leq k\leq s-1$, we see that
\begin{align}\label{eqn:76}
\pr\big(\overline\sbf(H_{n, s})=\bar\vt\big)
=
\frac{1}{\#\bH_{n, s}}
\sum\displaystyle_1\
\prod_{k=1}^s
\Bigg[
\Bigg(
\begin{array}{c}
\hskip-4pt z\big(\bar\vt; \ell_k\big)\hskip-7pt\phantom{.} \\
\hskip-4pt 2\hskip-7pt \phantom{.}
\end{array}
\Bigg)
+
\frac{1}{2}\cdot z\big(\bar\vt; \ell_k\big)\cdot\Big(z\big(\bar\vt; \ell_k-1\big)-1\Big)
\Bigg]
=:\frac{W_s(\bar\vt)}{\#\bH_{n, s}}\, .
\end{align}
Note that the quantity $W_s(\bar\vt)$ satisfies 
\[
0\leq W_1(\bar\vt)^s-s!\cdot W_s(\bar\vt)\leq 
c\cdot W_1(\bar\vt)^{s-1}\cdot\max_{\ell}z(\bar\vt; \ell)^2\, ,
\]
which essentially says that $n^{-3s/2}\cdot W_s(H_{n, 0})$ can be replaced by 
$n^{-3s/2}\cdot W_1(H_{n,0})^s/s!$ for obtaining distributional asymptotics.
Now the argument given below \eqref{eqn:65} for $s=1$ can be carried over to complete the proof for a general $s$.

\subsection{Proof of Theorem \ref{thm:bf=df-crum}}\label{sec:crum-proof}
Throughout this section we work with a fixed $g\geq 1$.
It is known \cite{goupil-schaeffer, chapuy} that 
\begin{gather}
\#\bum\sim\frac{4^g}{3^g\cdot g!\cdot\sqrt{\pi}}n^{3g-3/2}4^n\, .\label{eqn:93}
\end{gather}

Further, for any $s\geq 1$, $\bbf:\bM_{n, s}\to\bft(n, s)$ is a bijection, and consequently,  
\begin{align}\label{eqn:102}
\#\bM_{n, s}
&=\sum_{\vt\in\bM_{n, 0}}\#\bfac(\vt, s)
=\#\bM_{n, 0}\cdot\bE\big[\#\bfac(T_n, s)\big]
\sim\frac{\#\bM_{n, 0}}{s!}\cdot (2n)^{3s/2}\cdot\frac{\bE\big[B(C_n)^s\big]}{(2n)^{3s/2}}\\
&\sim
\frac{\#\bM_{n, 0}}{s!}\cdot (2n)^{3s/2}\cdot\bE\big[\big(\int_0^{\infty}\eta(\ve; 1, y)^2 dy\big)^s\big]
\sim\frac{2^{\frac{5s}{2}-2}}{s!\sqrt{\pi}}\cdot n^{\frac{3}{2}(s-1)}4^n\cdot\bE\big[\big(\int_0^1\ve(t)dt\big)^s\big]
\, ,\notag
\end{align}
where the third step uses Lemma \ref{lem:3}, the fourth step uses \eqref{eqn:29} and \eqref{eqn:44}, and the last step uses \eqref{eqn:74} and the fact that 
$\#\bM_{n, 0}=n^{-1}\cdot\dbinom{2n-2}{n-1}\sim 2^{2n-2}(n\pi)^{-1/2}$.

\begin{rem}
The enumeration technique in \eqref{eqn:102} is very much in the spirit of \cite{spencer-count}.
Let us point out here that this can also be deduced from the results of \cite{WaLe72}.
\greg{\cite[Formula  (5)]{WaLe72} states that the generating function $F_s(x)=\sum_{n\geq 0}\big(\#\bM_{n\sanch{+1},s}\big)\cdot x^n$ is of the form 
$F_0(x)(1-4x)^{-s}g_s(y)$. Here, $F_0(x)=(1-\sqrt{1-4x})/(2x)$ is the generating series for trees\footnote{This corrects a typo in \cite{WaLe72}, where the $x$ in the denominator of $F_0(x)$ is missing.}, $g_s$ is a polynomial of degree $(s-1)^+$ with positive integer coefficients, and $y=((1-4x)^{-1/2}-1)/2$. 
\sanch{(Recall that in our definition of $\bM_{n, s}$, we require the root vertex to have degree one. Thus, our notation is slightly different from that in \cite{WaLe72}, which explains why the coefficient of $x^n$ is $\#\bM_{n+1,s}$.)}
From this and elementary singularity analysis \cite[Theorem VI.1]{flajolet-sedgewick}, one obtains that as $n\to\infty$, for any $s\geq 1$ fixed, 
\begin{align}\label{eqn:1212}
\#\bM_{n,s} \sim \frac{\omega^*_s n^{\frac{3}{2}(s-1)}\sanch{4^n}}{\sanch{2^s}\cdot\Gamma((3s-1)/2)}\, ,
\end{align}
where $\omega^*_s$ is the leading coefficient of $g_s$. From  \cite[Formula (5a)]{WaLe72}, one immediately deduces that $\omega^*_s,s\geq 1$, satisfies the recursion
\begin{align}\label{eq:wright}
\omega^*_s=\sum_{k=1}^{s-1}\omega^*_k\omega^*_{s-k}+2(3s-4)\omega^*_{s-1}\, ,
\end{align}
which is exactly the recursion in \cite[Equation (13)]{janson2007brownian} (with the small difference that $\omega^*_0=1$ here rather  than $-1/2$ as in \cite{janson2007brownian}, but this does not matter since the initialization of \eqref{eq:wright} really starts at $s=1$). 
Now  \cite[Formula (14)]{janson2007brownian} identifies
\[\frac{\omega^*_s}{\Gamma((3s-1)/2)}=\frac{2^{\frac{5s}{2}-2}}{ s!\sqrt{\pi}}\cdot\bE\big[\big(\int_0^1 \sanch{2}\ve(t)dt\big)^s\big]\, ,
\]}
which combined with \eqref{eqn:1212} yields the asymptotics for $\#\bM_{n, s}$.
\end{rem}

Consider $\vm\in\bum$, and assume that the root vertex of $\vm$ has degree one.
Note also that $\spls(\vm)=2g$.
Thus, such a map $\vm\in\bum\cap\bM_{n, 2g}$.
We can explore $\vm$ in a breadth-first way as in Section \ref{sec:bf-df-search}.
Let $\bbf(\vm)=\big(\vt,\, (i_j, k_j;\, 1\leq j\leq 4g)\big)$.
Then the cardinality of the set of all such $\vm$ for which $i_1,\ldots, i_{4g}$ are not all distinct is $O(n^{-1/2})\cdot\#\bum$.
Indeed, this follows from the proof of Lemma \ref{lem:3} (which shows the analogous result for $\bM_{n, s}$) and the fact that $\#\bum=\Theta(\#\bM_{n, 2g})$ (which is a consequence of \eqref{eqn:93} and \eqref{eqn:102}). 
%(see \eqref{eqn:999} below for the expression for $\#\bum$).
Let $\bumast\subseteq\bum\cap\bM_{n, 2g}$ be the subset consisting of all $\vm$ for which $i_1,\ldots, i_{4g}$ are all distinct.
Removing the root vertex and the interior of the root edge of $\vm\in\bumast$ gives a map in $\bU\bM_{n-1, g}$.
Using this observation and \eqref{eqn:93}, it follows that
\begin{align}\label{eqn:95}
\#\bumast
\sim
\#\bU\bM_{n-1, g}
\sim
\frac{n^{3g-3/2}4^{g-1}}{3^g\cdot g!\cdot\sqrt{\pi}}4^n\, .
\end{align}
Clearly, $UM_{n, g}$ has the same scaling limit as a uniform element of $\bumast$, and so we can restrict our attention to unicellular maps in $\bumast$.

\begin{figure}
	\centering
	\includegraphics[trim=1cm 17cm 1cm 0cm, clip=true, angle=0, scale=.65]{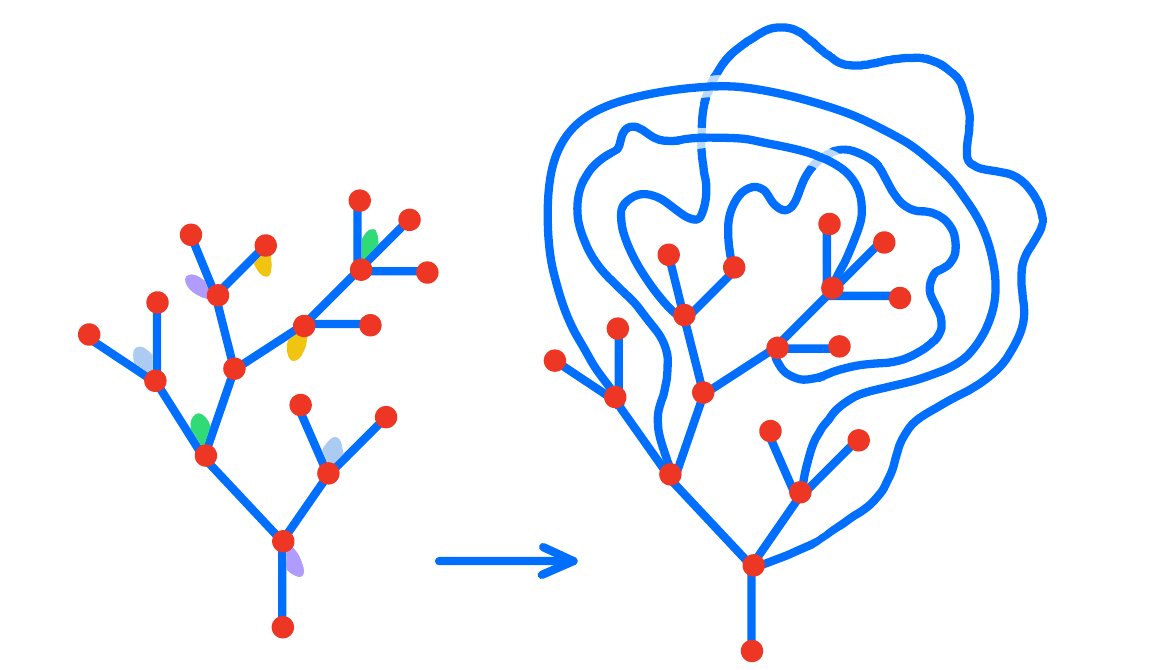}
	\captionof{figure}{Left: A plane tree with four pairs of corners corresponding to $\sigma=(1, 7)(2, 5)(3, 8)(4, 6)$ colored. Right: The associated unicellular map.}
	\label{fig:9}
\end{figure}
\begin{figure}
	\centering
	\includegraphics[trim=0cm 13cm 0cm 5cm, clip=true, angle=0, scale=.55]{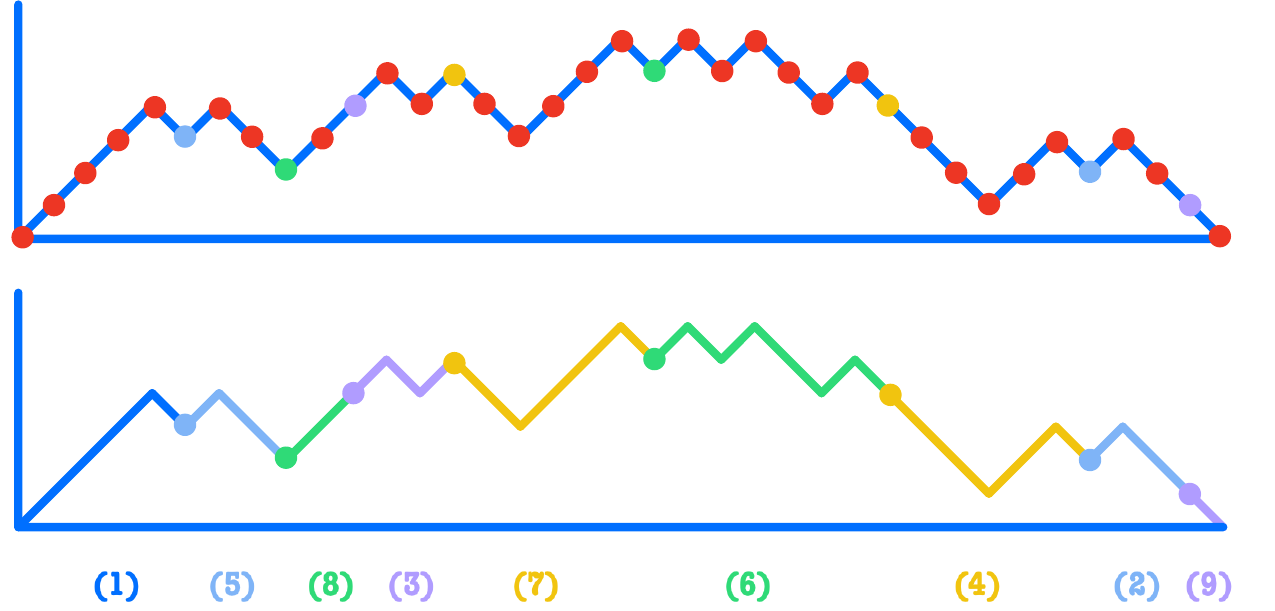}
	\captionof{figure}{Top: The contour function of the plane tree in Figure \ref{fig:9}. Bottom: Different segments of the same function numbered according to the order in which the face exploration of the unicellular map in Figure \ref{fig:9} traverses these parts.}
	\label{fig:10}
\end{figure}

Suppose $\sigma$ is a permutation on $[4g]$ whose cycle decomposition consists of $2g$ many transpositions: $\sigma=\prod_{j=1}^{2g}\big(\ell_{2j-1}, \ell_{2j}\big)$.
Suppose further that $\ell_1<\ell_3<\ldots<\ell_{4g-1}$, and $\ell_{2j-1}<\ell_{2j}$ for each $j\in [2g]$.
Recall the definition of the set $\bS_{(g)}$ from Section \ref{sec:def-S-g}.

\begin{lem}\label{lem:6}
	Suppose $\vt\in\bM_{n, 0}$, $\sigma$ is a permutation on $[4g]$ as above, and $ r_1< r_2<\ldots<r_{4g}$ are integers in $[2n-1]$.
	Let $i_j=r_{\ell_j}$ for $j\in[4g]$, and
	let $\vm\in\bM_{n, 2g}$ be the map obtained by adding an edge in $\vt$ between the $i_{2j-1}$-th corner and the $i_{2j}$-th corner of $\vt$, $1\leq j\leq 2g$, and endowing the neighbors of each vertex in $\vm$ with the natural circular order inherited from $\vt$.
	Then $\vm\in\bum\cap\bM_{n, 2g}$ iff $\sigma\in \bS_{(g)}$.
	
	Consequently, if $\sigma\in \bS_{(g)}$ and 
	$0\leq C_{\vt}\big(i_{2j-1}\big)-C_{\vt}\big(i_{2j}\big)\leq 1$ for $j\in [2g]$,
	then $\vm\in\bumast$ with $\bbf(\vm)=\big(\vt, (i_1,\ldots, i_{4g}, 1,\ldots, 1)\big)$.
\end{lem}

\chng{Before starting the proof, we ask the reader to recall the algebraic description of maps that was mentioned in Section~\ref{sec:notation}.
In particular, we will make use of the fact that a map is unicellular if and only if its face permutation consists of only one cycle.}

\vskip3pt

\noindent{\bf Proof of Lemma \ref{lem:6}:}
Let $e_1,\ldots, e_{2n}$ be the directed edges of $\vt$ in the contour order with $e_1$ being the root edge.
Then $\overrightarrow{E}(\vm)$ consists of $e_1,\ldots, e_{2n}$ together with an additional $4g$ directed edges; denote the directed edge going from the $r_j$-th corner to the $r_{\sigma(j)}$-th corner by $f_{j, \sigma(j)}$,
$j\in[4g]$.

Let us first assume that $\sigma\in \bS_{(g)}$.
As in Section \ref{sec:def-S-g}, let $\varrho$ be the permutation $\big(1, 2,\ldots, 4g\big)$.
Let $r_0=0$ and $r_{4g+1}=2n$.
Then the elements of $\overrightarrow{E}(\vm)$ in the contour exploration of $\vm$ starting with $e_1$ appear as follows:
For every $j\in\{0, 1, \ldots, 4g\}$, $e_{1+r_j}, e_{2+r_j}, \ldots, e_{r_{j+1}}$ appear consecutively in this order, 
and for every $j\in [4g]$, $e_{r_j}$ is followed by the directed edges $f_{j, \sigma(j)}$ and 
$e_{1+r_{\sigma(j)}}$.
(An illustration is given in Figures \ref{fig:9} and \ref{fig:10}.)
Using the fact that the permutation $\varrho\sigma$ has only one cycle of length $4g$, it follows that the face permutation of $\vm$ consists of only one cycle of length $2n+4g$, and consequently, $\vm$ is unicellular.

The converse follows by working the previous argument backwards.
If $\vm$ is unicellular, then its face permutation has only one cycle of length $2n+4g$.
Observing the order in which the directed edges of $\vm$ appear in its contour exploration, it follows that the permutation $\varrho\sigma$ must have only one cycle of length $4g$, and hence $\sigma\in \bS_{(g)}$. 

Now, because of the way $\ell_1,\ldots, \ell_{4g}$ are arranged among themselves, we have $i_{2j-1}<i_{2j}$ for $j\in[2g]$, and $i_1<i_3<\ldots<i_{4g-1}$.
Hence, if $0\leq C_{\vt}\big(i_{2j-1}\big)-C_{\vt}\big(i_{2j}\big)\leq 1$ for $j\in [2g]$,
then $ \big(i_1,\ldots, i_{4g}, 1,\ldots, 1\big)\in\bfac(\vt, 2g)$. 
Since $\cI\big(\vt, (i_1,\ldots, i_{4g}, 1,\ldots, 1)\big)=\vm$, it follows that $\bbf(\vm)=\big(\vt, (i_1,\ldots, i_{4g}, 1,\ldots, 1)\big)$.
If further $\sigma\in\bS_{(g)}$, then $\vm\in\bum\cap\bM_{n, 2g}$, and consequently $\vm\in\bumast$, as $i_1,\ldots, i_{4g}$ are distinct.
\qed

\medskip

For the rest of this section, the notation $\ell_1, \ldots, \ell_{2g}$ will be reserved for integers that are arranged among \chng{themselves} as described right before the statement of Lemma \ref{lem:6}, and  $\sigma=\prod_{j=1}^{2g}(\ell_{2j-1}, \ell_{2j})$ for a generic $\sigma\in \bS_{(g)}$.
%Further, we will use the notation $i_j=r_{\ell_j}$.
 
In view of Lemma \ref{lem:6}, $\bumast$ is in bijective correspondence with the set 
\begin{align}\label{eqn:96}
\big\{
\big(C, \sigma, (r_1,\ldots, r_{4g})\big)\, &:
\,\,
C\in\fC_n,\, \sigma\in \bS_{(g)},\, 1\leq r_1<\ldots r_{4g}\leq 2n-1\\
&\phantom{m}
\text{such that }
0\leq C\big(r_{\ell_{2j-1}}\big)-C\big(r_{\ell_{2j}}\big)\leq 1\text{ for }j\in[2g]
\big\}\, .\notag
\end{align}

Thus, sampling a uniform element of $\bumast$ is tantamount to sampling a uniform element in the above set, which can be done in the following steps:
For $C\in\fC_n$, $h\in\bZ_{>0}$, and $r, r'\in [2n-1]$, let 
$\psi_{C, h}(r, r')=\ind\big\{C(r)=h=C(r') \big\}+\ind\big\{C(r)=h, C(r')=h-1 \big\}$.

\vskip5pt

\begin{inparaenuma}
\noindent\item
Sample $C_n^{\dagger}$ with probability mass function (pmf) proportional to
\[
\sum_{\sigma\in \bS_{(g)}}\sum\displaystyle_2\sum\displaystyle_1
\prod_{j=1}^{2g}\psi_{C, h_j}\big(r_{\ell_{2j-1}},  r_{\ell_{2j}}\big)\, ,\ \ \ C\in\fC_n\, ,
\]
where $\sum_1$ and $\sum_2$ denote sum over all 
$1\leq r_1<r_2<\ldots <r_{4g}\leq 2n-1$ and $h_1,\ldots, h_{2g}\in\bZ_{>0}$ respectively.

\noindent\item 
Conditional on $C_n^{\dagger}=C$, sample $\Theta_n$ with pmf proportional to
\[
\sum\displaystyle_2\sum\displaystyle_1
\prod_{j=1}^{2g}\psi_{C, h_j}\big(r_{\ell_{2j-1}},  r_{\ell_{2j}}\big)\, ,\ \ \ \sigma\in \bS_{(g)}\, .
\]

\noindent\item 
Conditional on $C_n^{\dagger}=C$ and $\Theta_n=\sigma$, sample 
$\mvH_n=(H_{1, n},\ldots, H_{2g, n})$ with pmf proportional to
\[
\sum\displaystyle_1
\prod_{j=1}^{2g}\psi_{C, h_j}\big(r_{\ell_{2j-1}},  r_{\ell_{2j}}\big)\, ,\ \ \ h_1,\ldots, h_{2g}\in\bZ_{>0}\, .
\]

\noindent\item 
Conditional on $C_n^{\dagger}=C$, $\Theta_n=\sigma$, and $\mvH_n=(h_1,\ldots, h_{2g})$, sample 
$(v_{1, n},\ldots, v_{4g, n})$ with pmf proportional to
\[
\prod_{j=1}^{2g}\psi_{C, h_j}\big(r_{\ell_{2j-1}},  r_{\ell_{2j}}\big)\, ,\ \ \ 1\leq r_1<r_2<\ldots <r_{4g}\leq 2n-1\, .
\]
\end{inparaenuma}

\vskip5pt

Let $C_n$ be a uniform element of $\fC_n$.
Then using Proposition \ref{prop:local-time-convergence}, it follows that for $\sigma\in \bS_{(g)}$, as $n\to\infty$,
\begin{align}\label{eqn:94}
\frac{1}{(2n)^{3g}}\sum\displaystyle_2\sum\displaystyle_1
\prod_{j=1}^{2g}\psi_{C_n, h_j}\big(r_{\ell_{2j-1}},  r_{\ell_{2j}}\big)
\weakc
2^{2g}
\int_{(0,\infty)^{2g}}  \nu_{\ve, \sigma, \mvy}\big([0, 1]^{4g}\big) dy_1\ldots dy_{2g}
\, .
\end{align}
This explains the expression for the tilt in \eqref{eqn:78}.
Now, using the bijection explained around \eqref{eqn:96},
\begin{align}\label{eqn:97}
\#\bumast
=
\#\fC_n\cdot\bE\Big[
\sum_{\sigma\in\bS_{(g)}}\sum\displaystyle_2\sum\displaystyle_1
\prod_{j=1}^{2g}\psi_{C_n, h_j}\big(r_{\ell_{2j-1}},  r_{\ell_{2j}}\big)
\Big]\, .
\end{align}
Note that 
$\sum_2\sum_1\prod_{j=1}^{2g}\psi_{C_n, h_j}\big(r_{\ell_{2j-1}},  r_{\ell_{2j}}\big)
\leq
(2\|L_n\|_{\infty})^{2g}\cdot (2n)^{2g}$.
Thus, using \eqref{eqn:25}, it follows that we also have convergence of expectations in \eqref{eqn:94}.
Using this observation together with \eqref{eqn:97}, \eqref{eqn:95}, and the fact that 
$\#\fC_n=n^{-1}\cdot\dbinom{2n-2}{n-1}$, a direct calculation shows that
\begin{align}\label{eqn:98}
\sum_{\sigma\in\bS_{(g)}}
\bE\Big[\int_{(0,\infty)^{2g}}  \nu_{\ve, \sigma, \mvy}\big([0, 1]^{4g}\big) dy_1\ldots dy_{2g}\Big]
=
\Big(24^g\cdot g! \Big)^{-1}\, .
\end{align}
This explains the scaling constant in \eqref{eqn:78}.

Now using Proposition \ref{prop:local-time-convergence}  and arguments similar to the ones used in the proof of Proposition \ref{prop:bf-scaling-limit}, we can show that the following convergences happen jointly:
\begin{gather}
(2n)^{-1/2}C_n^{\dagger}(2n\,\cdot)\weakc\egum(\cdot)\, ,\ \ \ 
\Theta_n\weakc\Theta\, ,\ \ \ 
(2n)^{-1/2}\mvH_n\weakc\mvH\, ,\ \ \text{ and}\label{eqn:79}\\
(2n)^{-1}\big(v_{1, n},\ldots, v_{4g, n}\big)
\weakc
\big(u_1,\ldots, u_{4g}\big)\, ,\label{eqn:80}
\end{gather}
where the limiting random variables are as in Construction \ref{constr:crum-g-bf}.
The proof of Theorem \ref{thm:bf=df-crum} can now be completed using \eqref{eqn:79} and \eqref{eqn:80}.
We omit the details as no new idea is involved here.

\section{Discussion}\label{sec:discussion}
We will discuss some of the questions related to this work in this section.
Firstly, note that the result in Theorem \ref{thm:bf=df-Hs} concerns only the limiting space $\cH_{(s)}$.
However, our proof uses discrete approximation techniques.
It is natural to ask if this result can be proved directly in the continuum using properties of Brownian excursions.

As mentioned before, the proof of Theorem \ref{thm:bf=df-Hs} proceeds via a study of maps.
The reason is that an analogue of \eqref{eqn:24} for the contour function or the height function of uniform labeled trees, to the best of our knowledge, is not available in the literature.
More generally, one may hope for a result of the following form:

\begin{conj}
Let $\xi$, $T_{\xi}$, and $|T_{\xi}|$ be as in Theorem \ref{thm:conditioned-gw-tree}.
Write $\sigma^2=\var\xi$.
Let $H_n^{\xi}(t)$, $t=0, 1,\ldots, n$, (resp. $C_n^{\xi}(t)$, $t=0, 1, \ldots, 2n$) be the height function (resp. contour function) of $T_{\xi}$ conditioned to have $(n+1)$ vertices whenever $\pr\big(|T_{\xi}|=n+1\big)>0$.
Define a continuous function $L\big(H_n^{\xi}; t, y\big)$ on $[0, n]\times\bR$ 
(resp. a continuous function $L\big(C_n^{\xi}; t, y\big)$ on $[0, 2n]\times\bR$) 
by means of the formulas given in \eqref{eqn:21} and \eqref{eqn:22}.
Then 
\begin{gather*}
\Big(n^{-1/2}\cdot L\big(H_n^{\xi}; nt, y\sqrt{n}\big)\, ;\,  t\in[0, 1],\, y\in\bR\Big)
\weakc
\Big(\eta\Big(\frac{2\ve}{\sigma};\, t, y\Big)\, ;\,  t\in[0, 1],\, y\in\bR\Big)\, ,\text{ and}\\
\Big(n^{-1/2}\cdot L\big(C_n^{\xi}; 2nt, y\sqrt{n}\big)\, ;\,  t\in[0, 1],\, y\in\bR\Big)
\weakc
\Big(2\cdot\eta\Big(\frac{2\ve}{\sigma};\, t, y\Big)\, ;\,  t\in[0, 1],\, y\in\bR\Big)
\end{gather*}
in $C\big([0, 1]\times\bR\big)$, as $n\to\infty$ along the subsequence where $\pr\big(|T_{\xi}|=n+1\big)>0$.
\end{conj}

Such a result can be viewed as a generalization of \cite[Theorem 1.1]{drmota-gittenberger}.
A more ambitious project would be to identify the breadth-first construction of the stable graphs considered in \cite{GHS17}; 
see also \cite{SB-SD-vdH-SS, dhara-hofstad-leeuwaarden-sen-2}.
\cite[Theorem 3.2]{bhamidi-sen} gives an algorithm for constructing uniform connected graphs with a given degree sequence.
This algorithm can be thought of as a depth-first construction.
A similar algorithm can be developed from a breadth-first point of view.
One may try to use such an algorithm to obtain a breadth-first construction of the stable graphs studied in \cite{GHS17} (and thereby identifying its radius, two point function, and distance profile in terms of suitable functionals of a normalized excursion of an $\alpha$-stable L\'evy process). 
To carry out this program, one needs good control over the local time
field of the contour function (or height function) of the corresponding uniform plane tree with given (random) child sequence.
Although the question of the full breadth-first construction in the stable setting remains open, some partial results in this direction were obtained in the recent work \cite{clancy2021epidemics}.

Explicit expressions are known for the densities of the two point function and the radius of the space $\cT_{\ve}$.
The two point function of $\cT_{\ve}$ follows the Rayleigh distribution \cite{aldous-crt-1, aldous-crt-3}, whereas $\radi(\cT_{\ve})$ follows the more complicated Theta distribution \cite{biane-pitman-yor}\cite[Chapter V.4.3]{flajolet-sedgewick}.
In \cite{wang-minmin-height}, an expression (given in terms of an infinite series) for the joint distribution function of $\radi(\cT_{\ve})$ and the diameter of $\cT_{\ve}$ is computed using probabilistic arguments.
It would be interesting to see if such explicit expressions can be obtained for the laws of the radius and the two point function of the spaces $\cH_{(s)}$.

\appendix
\section{} \label{sec:appendix}
Our aim in this section is to outline a proof of \eqref{eqn:24}.
For $n\geq 1$, let $\big(S_n(j),\, 0\leq j\leq 2n+1\big)$ be a simple symmetric random walk of length $2n+1$ started at the origin, i.e., $S_n(0)=0$ and for $0\leq j\leq 2n$,
\[
\pr\big(S_n(j+1)=y\ \big|\ S_n(k),\, 0\leq k\leq j\big)=1/2\ \text{ if }\ y=S_n(j)\pm 1\, .
\]
Let 
\[
S_n^{\br}=\big(S_n\ \big|\ S_n(2n+1)=-1\big)\ \ \text{ and }\ \ 
S_n^{\exx}=\big(S_n^{\br}\ \big|\ S_n^{\br}(j)\geq 0\,\text{ for }\, 0\leq j\leq 2n\big)
\]
be the corresponding bridge and excursion of length $2n+1$ respectively.
We extend $S_n$, $S_n^{\br}$, and $S_n^{\exx}$ to continuous functions on $[0, 2n+1]$ by linear interpolation.
For $t\in[0, 2n+1]$ and $y\in\bR$, define $L(S_n; t, y)$ similar to \eqref{eqn:21} and \eqref{eqn:22}. 
Let
\begin{gather*}
\ell_n(\cdot, \cdot):=L(S_n; \cdot, \cdot)\, ,\ \ \ 
\bar S_n(t):=(2n+1)^{-1/2}\cdot S_n\big((2n+1)t\big)\, ,\ 0\leq t\leq 1
\, , \text{ and}\\
\bar\ell_n(t, y):=(2n+1)^{-1/2}\cdot\ell_n\big((2n+1)t,\, y\sqrt{2n+1}\big)\, ,\ \ \ 
0\leq t\leq 1\,\ y\in\bR\, .
\end{gather*}
Similarly define 
$\ell_n^{\br}$, $\bar S_n^{\br}$, $\bar\ell_n^{\br}$, $\ell_n^{\exx}$, $\bar S_n^{\exx}$, and $\bar\ell_n^{\exx}$.
Note that the function $\big(C_n(t+1)-1\, ,\ 0\leq t\leq 2n-1\big)$ has the same law as $S_{n-1}^{\exx}$.
Consequently, \eqref{eqn:24} will follow if we show that
\begin{align}\label{eqn:81}
\big(\bar S_n^{\exx}\, ,\ \bar\ell_n^{\exx} \big)
\weakc
\big(\ve\, ,\, \eta(\ve;\cdot, \cdot)\big)
\end{align}
in $C[0, 1]\times C\big([0, 1]\times\bR\big)$.

Let $\big(B(t)\, ,\ 0\leq t\leq 1\big)$ (resp. $\big(B^{\br}(t)\, ,\ 0\leq t\leq 1 \big)$) be a standard one dimensional Brownian motion started at the origin (resp. Brownian bridge with $B^{\br}(0)=B^{\br}(1)=0$).
For the rest of the proof, $\eta^{\bm}$ and $\eta^{\br}$ will respectively denote 
$\eta(B; \cdot, \cdot)$ and $\eta(B^{\br}; \cdot, \cdot)$.

By \cite[Theorem F]{csaki-revesz} (see also \cite{revesz, bass-khoshnevisan}), it follows that we can construct a standard one dimensional Brownian motion $W(t)$, $t\geq 0$, and $S_n$ for all $n\geq 1$ simultaneously on a single probability space such that for any $\delta>0$,
\begin{align}\label{eqn:82}
n^{-1/4-\delta}
\Big(
\max_{0\leq j\leq 2n+1}\max_{y\in \bZ}
\big|\ell_n(j, y)-\eta(W; j, y)\big|
+
\max_{0\leq j\leq 2n+1}
\big|S_n(j)-W(j)\big|
\Big)
\convas 0\, ,
\end{align}
as $n\to\infty$.
Writing $\max_1$ for maximum over 
$t\in\big\{1/(2n+1), 2/(2n+1),\ldots, 1 \big\}$ and $y\in\bZ/\sqrt{2n+1}$, it follows from \eqref{eqn:82} that in this space, for any $\delta\in(0, 1/4)$,
\begin{align}\label{eqn:83}
&\max\displaystyle_1
\Big|
\bar \ell_n(t, y)-\eta\big(\bar W; t, y\big)
\Big|\notag\\
&\hskip40pt
=
\max\displaystyle_1
\frac{1}{\sqrt{2n+1}}\Big|
\ell_n\big((2n+1)t, y\sqrt{2n+1}\big)-\eta\Big(W; (2n+1)t, y\sqrt{2n+1}\Big)
\Big|\notag\\
&\hskip80pt
=
O(n^{-1/4+\delta})\, ,\ \text{ a.s.}\, ,
\end{align}
where $\bar W(t)=(2n+1)^{-1/2}W\big((2n+1)t\big)$, $0\leq t\leq 1$.
Further, in this space,
\begin{align}\label{eqn:84}
\max_{t\in\big\{1/(2n+1),\ 2/(2n+1),\ldots,\ 1\big\}}
\big|\bar S_n(t)-\bar W(t)\big|=
O(n^{-1/4+\delta})\, ,\ \text{ a.s.}\, .
\end{align}
Since $\bar W\equald B$, \eqref{eqn:83} and \eqref{eqn:84} yield 
\begin{align}\label{eqn:85}
\big(\bar S_n, \bar\ell_n\big)
\weakc
\big(B, \eta^{\bm}\big)
\end{align}
in $C[0, 1]\times C\big([0, 1]\times\bR\big)$.

Fix $\eps\in(0, 1)$, and let $B_{\eps}$ and $\eta_{\eps}^{\bm}$ denote the restrictions of $B$ and $\eta^{\bm}$ to $[0, 1-\eps]$ and $[0, 1-\eps]\times\bR$ respectively.
Similarly define $B_{\eps}^{\br}$ and $\eta_{\eps}^{\br}$.
Then
\begin{align}\label{eqn:86}
\bE\big[\phi\big(B_{\eps}^{\br}, \eta_{\eps}^{\br}\big)\big]
=
\frac{
\bE\big[\phi\big(B_{\eps}, \eta_{\eps}^{\bm}\big)\exp\big(-B(1-\eps)^2/(2\eps)\big)\big]
}{
\bE\big[\exp\big(-B(1-\eps)^2/(2\eps)\big)\big]
}
\end{align}
for any bounded continuous $\phi$.

Let $2m+1$ be the smallest odd integer bigger than $(1-\eps)(2n+1)$. 
Let $S_{n,m}^{\br}$ and $\ell_{n, m}^{\br}$ denote the restrictions of $S_n^{\br}$ and $\ell_n^{\br}$ to $[0, 2m+1]$ and $[0, 2m+1]\times\bR$ respectively.
Then
\begin{align}\label{eqn:87}
\bE\big[\phi\big(S_{n,m}^{\br}, \ell_{n,m}^{\br}\big)\big]
=
\frac{
	\bE\big[\phi\big(S_m, \ell_m\big)\pr\big(X(2n-2m)=S_m(2m+1)+1\ \big|\ S_m\big)\big]
}{
	\bE\big[\pr\big(X(2n-2m)=S_m(2m+1)+1\ \big|\ S_m\big)\big]
}\, ,
\end{align}
where $X(2n-2m)$ is a simple symmetric random walk started at the origin and run up to time $2n-2m$ independent of all other random variables.
By the local central limit theorem,
\begin{align}\label{eqn:88}
\sup_{j\in\bZ}
\Big|\pr\big(X(2n-2m)=2j\big)-\big(4\pi(n-m)\big)^{-1/2}\exp\Big(-\frac{4j^2}{4(n-m)}\Big)\Big|
=
o(n^{-1/2})\, .
\end{align}
It follows from \eqref{eqn:85}, \eqref{eqn:86}, \eqref{eqn:87}, and \eqref{eqn:88} that for every $\eps\in(0, 1)$,
\begin{align}\label{eqn:89}
\big(\bar S_{n, \eps}^{\br},\ \bar\ell_{n, \eps}^{\br}\big)
\weakc
\big(B_{\eps}^{\br}, \eta_{\eps}^{\br}\big)\, ,
\end{align}
where $\bar S_{n, \eps}^{\br}$ and $\bar\ell_{n, \eps}^{\br}$ are restrictions of 
$\bar S_n^{\br}$ and $\bar\ell_n^{\br}$ to $[0, 1-\eps]$ and $[0, 1-\eps]\times\bR$ respectively.
It would follow from \eqref{eqn:89} that
\begin{align}\label{eqn:90}
\big(\bar S_n^{\br},\ \bar\ell_n^{\br}\big)
\weakc
\big(B^{\br}, \eta^{\br}\big)\, ,
\end{align}
provided we could show that for every $\delta_1, \delta_2>0$ there exists $\eps_0>0$ such that for $\eps\in(0, \eps_0)$,
\begin{align}\label{eqn:91}
\limsup_{n\to\infty}\
\pr\Big(\sup_{1-\eps\leq t\leq 1}|\bar S_n^{\br}(t)|
+
\sup_{y\in\bR}\big(\bar\ell_n^{\br}(1, y)-\bar\ell_n^{\br}(1-\eps, y)\big)>\delta_1\Big)
\leq\delta_2\, .
\end{align}
However, \eqref{eqn:91} is immediate upon observing that the time reversal of the process $-(1+S_n^{\br})$ has the same law as $S_n^{\br}$ and then using \eqref{eqn:89}.

Now the Vervaat transform of $B^{\br}$ (resp. $S_n^{\br}$) with respect to its almost surely unique global minima (resp. the first global minima) has the same law as $\ve$ (resp. $S_n^{\exx}$).
Using this fact together with \eqref{eqn:90}, \eqref{eqn:81} follows.

\section*{Acknowledgments}
The authors thank Louigi Addario-Berry for many insightful discussions about the results of \cite{lab-oa-gc-ef-cg}.
The authors also thank two anonymous referees for their careful reading and many comments and suggestions on an earlier version of the paper.
Part of the work was done when both authors were attending Oberwolfach workshop 1750-Network Models: Structure and Function.
The authors thank the organizers of the workshop and the Oberwolfach Research Institute for Mathematics for a stimulating work environment. GM acknowledges support from Institut Universitaire de France and Fondation Simone et Cino Del Duca. 
SS was partially supported by MATRICS grant MTR/2019/000745 from SERB, and by the Infosys Foundation, Bangalore.

\bibliographystyle{plain}
\bibliography{breadth_first_bib}

\end{document}